%% file: apos.tex
\documentstyle[a4,txmac,%
case,%
mypic,%
epsf,%
twoside,%
nocaphead,%
amssymb,times,mathptm]{article}

\advance\oddsidemargin by -0.7cm
\advance\evensidemargin by -0.7cm
\advance\textwidth by 1.4cm

\def\mynewtheo#1#2{%
\newtheorem{@#1}{#2}[section]%
\newenvironment{#1}{\begin{@#1}\rm}{\end{@#1}}}
 
\mynewtheo{lem}{Lemma}
\mynewtheo{exer}{Exercise}
\mynewtheo{theo}{Theorem}
\mynewtheo{rem}{Remark}
\mynewtheo{defi}{Definition}
\mynewtheo{conj}{Conjecture}
\mynewtheo{corr}{Corollary}
\mynewtheo{prop}{Proposition}
\mynewtheo{ques}{Question}
\mynewtheo{exam}{Example}

\newenvironment{eqn}{\begin{equation}}{\end{equation}}
\newenvironment{theorem}{\begin{theo}}{\end{theo}}

\def\eqref#1{\mbox{(\protect\reference{#1}})}
\def\proof{\@ifnextchar[{\@proof}{\@proof[\unskip]}}
\def\@proof[#1]{\noindent{\bf Proof #1.}\enspace}

\begin{document}

\makeatletter

\parskip5pt plus 1pt minus 2pt

\author{A. Stoimenow\\[2mm]
\small Humboldt University Berlin, Dept.~of Mathematics,\\
\small Ziegelstra\ss e 13a, 10099 Berlin, Germany,\\
\small e-mail: {\tt stoimeno@informatik.hu-berlin.de},\\
\small WWW: {\hbox{\tt http://www.informatik.hu-berlin.de/\raisebox{-0.8ex}{\tt\~{}}stoimeno}}
}

\title{\large\bf
\uppercase{Gauss sums on almost positive knots}}

\date{\large Current version: \today\ \ \ First version:
\makedate{23}{10}{1997}}

\maketitle

\def\chrd#1#2{\picline{1 #1 polar}{1 #2 polar}}
\def\arrow#1#2{\picvecline{1 #1 polar}{1 #2 polar}}

\def\labch#1#2#3{\chrd{#1}{#2}\picputtext{1.3 #2 polar}{$#3$}}
\def\labar#1#2#3{\arrow{#1}{#2}\picputtext{1.3 #2 polar}{$#3$}}
\def\labbr#1#2#3{\arrow{#1}{#2}\picputtext{1.3 #1 polar}{$#3$}}

\def\CD{\szCD{6mm}}
\def\szCD#1#2{{\let\@nomath\@gobble\small\diag{#1}{2.4}{2.4}{
  \pictranslate{1.2 1.2}{
    \piccircle{0 0}{1}{}
    #2
}}}}
\let\GD\CD

\def\@dcont{}
\def\svCD#1{\ea\glet\csname #1\endcsname\@dcont}
\def\rsCD#1{\ea\glet\ea\@dcont\csname #1\endcsname\ea\glet
\csname #1\endcsname\relax}

\def\addCD#1{\ea\gdef\ea\@dcont\ea{\@dcont #1}}
\def\drawCD#1{\szCD{#1}{\@dcont}}

%

\long\def\@makecaption#1#2{%
   \vskip 10pt
   {\let\label\@gobble
   \let\ignorespaces\@empty
   \xdef\@tempt{#2}%
   }%
   \ea\@ifempty\ea{\@tempt}{%
   \setbox\@tempboxa\hbox{%
      \fignr#1#2}%
      }{%
   \setbox\@tempboxa\hbox{%
      {\fignr#1:}\capt\ #2}%
      }%
   \ifdim \wd\@tempboxa >\captionwidth {%
      \rightskip=\@captionmargin\leftskip=\@captionmargin
      \unhbox\@tempboxa\par}%
   \else
      \hbox to\captionwidth{\hfil\box\@tempboxa\hfil}%
   \fi}%
\def\fignr{\small\sffamily\bfseries}%
\def\capt{\small\sffamily}%

\newdimen\@captionmargin\@captionmargin2\parindent\relax
\newdimen\captionwidth\captionwidth\hsize\relax

\def\nin{\not\in}
\def\bC{{\Bbb C}}
\def\bR{{\Bbb R}}
\def\bN{{\Bbb N}}
\def\cK{{\cal K}}
\def\bt{\bar t_2}

\def\pr{\text{\rm pr}\,}
\def\ncap{\not\mathrel{\cap}}
\def\|{\mathrel{\kern1.5pt\Vert\kern1.5pt}}
\def\lra{\longrightarrow}
\let\ds\displaystyle
\let\reference\ref
\def\so{\Rightarrow}
\let\ex\exists
\let\fa\forall

\def\bysame{\same[\kern2cm]\,}
\def\qed{\hfill\@mt{\Box}}
\def\@mt#1{\ifmmode#1\else$#1$\fi}

\def\proof{\@ifnextchar[{\@proof}{\@proof[\unskip]}}
\def\@proof[#1]{\noindent{\bf Proof #1.}\enspace}

\def\myfrac#1#2{\raisebox{0.2em}{\small$#1$}\!/\!\raisebox{-0.2em}{\small$#2$}}
\newcommand{\mybr}[2]{\text{$\Bigl\lfloor\mbox{%
\small$\displaystyle\frac{#1}{#2}$}\Bigr\rfloor$}}
\def\mybrtwo#1{\mbox{\mybr{#1}{2}}}

\let\x\exists
\let\fa\forall
\let\sg\sigma
\let\tl\tilde
\let\dl\delta
\let\Dl\Delta
\let\eps\varepsilon
\let\ol\overline
\def\bt{\bar t_2}

\def\epsfs#1#2{{\ifautoepsf\unitxsize#1\relax\else
\epsfxsize#1\relax\fi\epsffile{#2.eps}}}
\def\epsfsv#1#2{{\vcbox{\epsfs{#1}{#2}}}}
\def\vcbox#1{\setbox\@tempboxa=\hbox{#1}\parbox{\wd\@tempboxa}{\box
  \@tempboxa}}

\def\@test#1#2#3#4{%
  \let\@tempa\go@
  \@tempdima#1\relax\@tempdimb#3\@tempdima\relax\@tempdima#4\unitxsize\relax
  \ifdim \@tempdimb>\z@\relax
    \ifdim \@tempdimb<#2%
      \def\@tempa{\@test{#1}{#2}}%
    \fi
  \fi
  \@tempa
}

\def\go@#1\@end{}
\newdimen\unitxsize
\newif\ifautoepsf\autoepsftrue

\unitxsize4cm\relax
\def\epsfsize#1#2{\epsfxsize\relax\ifautoepsf
  {\@test{#1}{#2}{0.1 }{4   }
		{0.2 }{3   }
		{0.3 }{2   }
		{0.4 }{1.7 }
		{0.5 }{1.5 }
		{0.6 }{1.4 }
		{0.7 }{1.3 }
		{0.8 }{1.2 }
		{0.9 }{1.1 }
		{1.1 }{1.  }
		{1.2 }{0.9 }
		{1.4 }{0.8 }
		{1.6 }{0.75}
		{2.  }{0.7 }
		{2.25}{0.6 }
		{3   }{0.55}
		{5   }{0.5 }
		{10  }{0.33}
		{-1  }{0.25}\@end
		\ea}\ea\epsfxsize\the\@tempdima\relax
		\fi
		}

\input{myeqn.tex}

\let\diagram\diag

\def\boxed#1{\diagram{1em}{1}{1}{\picbox{0.5 0.5}{1.0 1.0}{#1}}}

\def\rato#1{\hbox to #1{\rightarrowfill}}
\def\arrowname#1{{\enspace
\setbox7=\hbox{F}\setbox6=\hbox{%
\setbox0=\hbox{\footnotesize $#1$}\setbox1=\hbox{$\to$}%
\dimen@\wd0\advance\dimen@ by 0.66\wd1\relax
$\stackrel{\rato{\dimen@}}{\copy0}$}%
\ifdim\ht6>\ht7\dimen@\ht7\advance\dimen@ by -\ht6\else
\dimen@\z@\fi\raise\dimen@\box6\enspace}}

\def\contr{\diagram{1em}{0.6}{1}{\piclinewidth{35}%
\picstroke{\picline{0.5 1}{0.2 0.4}%
\piclineto{0.6 0.6}\picveclineto{0.3 0}}}}

\def\abstractname{}

{\let\@noitemerr\relax
\vskip-2.7em\kern0pt\begin{abstract}
\noindent{\bf Abstract.}\enspace
Using the Fiedler-Polyak-Viro Gau\ss{} diagram formulas we study
the Vassiliev invariants of degree $2$ and $3$ on almost positive knots.
As a consequence we show that the number of almost positive knots
of given genus or unknotting number
grows polynomially in the crossing number, and also recover and extend,
\em{inter alia} to their untwisted Whitehead doubles, previous results 
on the polynomials and signatures of such knots. In particular,
we prove that there are no achiral almost positive knots and
classify all almost positive diagrams of the unknot.
We give an application to contact geometry (Legendrian knots)
and property $P$.
\end{abstract}
}
\vspace{1mm}

\section{Introduction}
  
Many properties of knots are defined by the
existence of diagrams with certain properties. Such classical properties
are alternation and positivity. Adjoining the word ``almost''
before the name of the property, we mean that the knot does
not have a diagram with that property, but one in which it can
be attained by one crossing change.
Here we consider the notion for positivity.

\begin{defi}
The writhe is a number ($\pm1$), assigned to any crossing in a link
diagram. A crossing as in figure \ref{figwr}(a) has writhe 1 and
is called positive. A crossing as in figure \ref{figwr}(b) has writhe 
$-1$ and is called negative.

\begin{figure}[htb]
\[
\begin{array}{c@{\qquad}c}
\diag{6mm}{1}{1}{
\picmultivecline{0.12 1 -1.0 0}{1 0}{0 1}
\picmultivecline{0.12 1 -1.0 0}{0 0}{1 1}
} &
\diag{6mm}{1}{1}{
\picmultivecline{0.12 1 -1 0}{0 0}{1 1}
\picmultivecline{0.12 1 -1 0}{1 0}{0 1}
}
\\[2mm]
(a) & (b)
\end{array}
\]
\caption{\label{figwr}}
\end{figure}
\end{defi}

\begin{defi}
A knot is called positive if it has a positive diagram,
i.~e. a diagram with all crossings positive.
A knot is called almost positive if it is not positive, but
has an almost positive diagram, i.~e. a diagram with all crossings
positive except one.
\end{defi}


Recently, the Fiedler-Polyak-Viro approach \cite{VirPol,Fied2}
to Vassiliev invariants \cite{BarNatanVI,BarNatanBibl,Vassiliev}
via Gau\ss{} diagram formulas gave a new powerful tool in
studying positivity, see \cite{pos}.

The aim of the present paper is to extend the applications of
the Fiedler Gau\ss{} sum formula to almost positivity. In particular we
will apply this formula to classify almost positive diagrams of the
unknot (corollary \reference{Cr3.2}). More generally, we will show that
knots with zero or negative Fiedler invariant 
cannot be almost positive. As achiral knots have zero Fiedler invariant,
this means in particular that any almost positive knot is chiral.

We also show chirality for the untwisted Whitehead doubles of an
almost positive knot (with either clasps), which follows from the
positivity of its Casson invariant $v_2$ (\S\reference{SCas}). 
%
As a further consequence we prove that an almost positive knot itself
has non-trivial polynomial invariants. This is proved for the 
Jones polynomial (and hence also for the HOMFLY \cite{HOMFLY}
and Kauffman \cite{Kauffman} $F$ polynomial), but also for the
Alexander polynomial \cite{Alexander} and the $Q$ polynomial of
Brandt-Lickorish-Millett \cite{BLM} and Ho \cite{Ho}.

In \S\reference{Sgen} we improve the positivity result for $v_2$
on almost positive knots of given genus to an estimate involving
the crossing number, similar to the one on positive knots in
\cite{pos}. This allows to extend a result of \cite{gen1} to
almost positive knots.

\begin{theorem}\label{tgu}
The number of almost positive knots of given genus of unknotting number
grows polynomially in the crossing number. That is,
\[
\#\,\{\,K\,:\,\mbox{$K$ almost positive knot},\,c(K)=n,\,g(K)=g\,\}\,=\,
O_n(n^{p_g})\,,
\]
for some number $p_g\in\bN$, where $O_n$ denotes the asymptotic
behaviour as $n\to\infty$ (and the same statement holds for `$g(K)$'
replaced by `$u(K)$').
\end{theorem}

The positivity result for $v_2$ has a consequence to the behaviour of
the degrees of the skein polynomial of almost positive knots. An
application of this is the generalization (to almost positive knots)
of the result of Kanda \cite{Kanda} (see also \cite{CGM,TabFuchs}),
which followed in more general form for positive knots already from
\cite{beha}. (Here, and in the sequel, for a knot $K$, 
$!K$ denotes the obverse, or mirror image, of $K$.
The definition of $tb$ and $\mu$ is recalled in \S\reference{Sgen}.)

\begin{theorem}\label{tL}
Let $\{K_i\}$ be distinct positive or almost positive knots and
$\{\cK_i\}$ be Legendrian embeddings of $!K_i$ in the 
in the standard contact space $(\bR^3(x,y,z),\,dx+y\,dz)$ with
Thurston-Bennequin invariants $tb(\cK_i)$ and
Maslov (rotation) numbers $\mu(\cK_i)$. Then 
\[
tb(\cK_i)+|\mu(\cK_i)|\,\arrowname{i\to\infty}\,-\infty\,.
\]
The same statement holds for $\{\cK_i\}$ transverse, when omitting
the $\mu(\cK_i)$ term.
\end{theorem}

Finally, in \S\reference{Ssig} we give a proof that almost positive
knots have positive signature.

\section{Gau\ss{} sums}

The ``almost'' concept was introduced by C.\ Adams et al.\ \cite{Adams2}
for alternation. Almost alternating knots are much more common than
almost positive, but their diversity makes proofs of specific
properties difficult. 

\begin{figure}[htb]
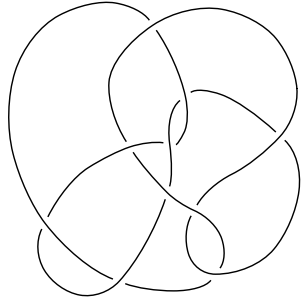

\[
\epsfs{4cm}{k-10-145a}
\]
\caption{\label{fig10-145}The knot $!10_{145}$ is almost positive.}
\end{figure}

The simplest example of an almost positive knot is given in
figure \reference{fig10-145}.

\begin{exam}
The knot $!10_{145}$ of \cite{Rolfsen} is almost positive, as
its diagram on figure \reference{fig10-145} is so, but it is known
not to be positive \cite{Cromwell}.
\end{exam}

We use the Alexander-Briggs notation and the Rolfsen \cite{Rolfsen}
tables to distinguish between a knot and its obverse. ``Projection''
is the same as ``diagram'', and this means a knot or link
diagram. Diagrams are always assumed oriented.

So far much less than on almost alternating knots seems to be known on
almost positive knots. However, the concept of Gau\ss{} sum invariants
developed by Fiedler-Polyak-Viro \cite{Fied2,VirPol} has
several direct applications to such knots. We recall the basic
parts of this concept now.

\begin{@defi}[\protect\cite{Fied2}]\rm
A Gau\ss{} diagram of a knot diagram is an oriented circle with
arrows connecting points on it mapped to a crossing and
oriented from the preimage of the undercrossing to the
preimage of the overcrossing.
\end{@defi}

\begin{figure}[htb]
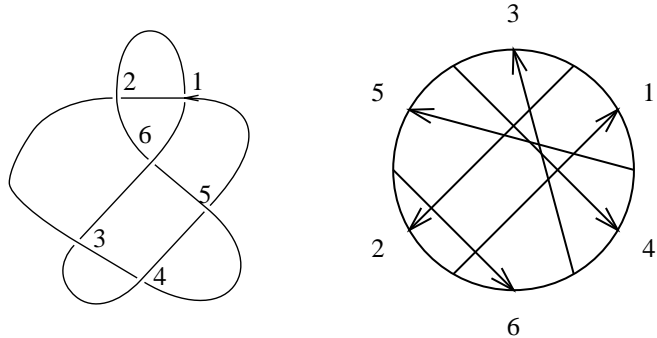

{
\[
\begin{array}{c@{\qquad}c}
     \diag{2cm}{2.0}{2.0}{
  \picputtext[dl]{0 0}{\autoepsffalse\epsfs{4cm}{6_2}}
  \picputtext[u]{1.25 1.65}{$1$} 
  \picputtext[u]{0.8 1.65}{$2$} 
  \picputtext[u]{0.6 0.6}{$3$} 
  \picputtext[u]{1.0 0.35}{$4$} 
  \picputtext[u]{1.3 0.9}{$5$} 
  \picputtext[u]{0.9 1.3}{$6$} 
	}
 & 
\diag{1.6cm}{2.4}{2.4}{\pictranslate{1.2 1.2}{
\piclinewidth{60}
\piccircle{0 0}{1}{}
\labar{240}{ 30}{1}
\labar{ 60}{210}{2}
\labar{300}{ 90}{3}
\labar{120}{330}{4}
\labar{  0}{150}{5}
\labar{180}{270}{6}
} 
} 
\end{array}
\]}
\caption{The standard diagram of the knot $6_2$ and its Gau\ss{}
diagram.\label{fig6_2}}
\end{figure}

\begin{exam}
As an example, figure \ref{fig6_2} shows the knot $6_2$ in its
commonly known projection and the corresponding Gau\ss{} diagram.
\end{exam}

Fiedler \cite{Fied2,NewInv} found the following formula for (a
variation of) the degree-3-Vassiliev invariant using Gau\ss{} sums.
\begin{eqn}\label{v3}
v_3\,=\,\sum_{(3,3)}w_pw_qw_r\,+\,
\sum_{(4,2)0}w_pw_qw_r\,+\,\frac{1}{2}\sum_\text{$p,q$ linked}(w_p+w_q)\,,
\end{eqn}
where the configurations are
\[
\begin{array}{c@{\qquad}c@{\qquad}c}
\CD{\chrd{50}{240}
\chrd{90}{270}
\chrd{125}{320}} &
\CD{\chrd{-180}{0}
\arrow{-70}{70}
\arrow{-110}{110}} &
\CD{\arrow{250}{70}
\arrow{290}{110}} \\[4mm]
(3,3) & (4,2)0 & \text{$p,q$ linked}
\end{array}
\]
Here chords depict arrows which may point in both directions and
$w_p$ denotes the writhe of the crossing $p$. For a given configuration,
the summation in \eqref{v3} is done over each unordered pair/triple of
crossings, whose arrows in the Gau\ss{} diagram form that configuration.
If two chords intersect, we call the corresponding
crossings \em{linked}, and the crossing, whose over-crossing of followed
by the under-crossing of the other when passing the diagram in the
orientation direction we call \em{distinguished}. We will denote
by $lk(D)$ the number of linked pairs in a diagram $D$.

To make precise which variation of the degree-3-Vassiliev invariant
we mean, we noted in \cite{NewInv}, that
\[
v_3\,=\,-\frac 13V^{(2)}(1)-\frac 19V^{(3)}(1)\,,
\]
where $V$ is the Jones polynomial \cite{Jones}. We noted also,
that $v_3$ is asymmetric, i.~e., $v_3(!K)=-v_3(K)$, so that achiral
knots have zero invariant.

\begin{defi}
The diagram on the left of figure \reference{figtan}
is called \em{connected sum} $A\# B$ of the diagrams $A$ and $B$.
If a diagram $D$ can be represented as the connected sum of 
diagrams $A$ and $B$, such that both $A$ and $B$ have at least one
crossing, then $D$ is called \em{disconnected} (or composite), else
it is called \em{connected} (or prime).
\end{defi}

\begin{figure}[htb]
\[
\diag{6mm}{3}{2}{
  \piccirclearc{1.8 1}{0.5}{-120 120}
  \picfilledcircle{1 1}{0.8}{$A$}
}\, \#\,
\diag{6mm}{3}{2}{
  \piccirclearc{1.2 1}{0.5}{60 300}
  \picfilledcircle{2 1}{0.8}{$B$}
}\quad =\quad
\diag{6mm}{4}{2}{
  \piccirclearc{2 0.5}{1.3}{45 135}
  \piccirclearc{2 1.5}{1.3}{-135 -45}
  \picfilledcircle{1 1}{0.8}{$A$}
  \picfilledcircle{3 1}{0.8}{$B$}
} 
\]
\caption{\label{figtan}}
\end{figure}

\begin{defi}
A crossing $p$ in a knot diagram $D$ is called \em{reducible}
(or nugatory) if it is linked with no other crossing. Then $D$
looks like on the left of figure \reference{figred}. $D$
is called reducible if it has a reducible crossing, else it is
called \em{reduced}. The reducing of the reducible crossing $p$
is the move depicted on figure \reference{figred}. Each diagram $D$
can be (made) reduced by a finite number of these moves.
\end{defi}

\begin{figure}[htb]
\begin{eqn}\label{eqred}
\diag{6mm}{4}{2}{
  \picrotate{-90}{\rbraid{-1 2}{1 1.4}}
  \picputtext[u]{2 0.7}{$p$}
  \picscale{1 -1}{
    \picfilledcircle{0.7 -1}{0.8}{$P$}
  }
  \picfilledcircle{3.3 1}{0.8}{$Q$}
} \qquad\lra\qquad
\diag{6mm}{4}{2}{
  \piccirclearc{2 0.5}{1.3}{45 135}
  \piccirclearc{2 1.5}{1.3}{-135 -45}
  \picfilledcircle{1 1}{0.8}{$P$}
  \picfilledcircle{3 1}{0.8}{$Q$}
} 
\end{eqn}
\caption{\label{figred}}
\end{figure}

\begin{defi}
For two chords in a Gau\ss{} diagram $a\cap b$ means
``$a$ intersects $b$'' (or crossings $a$ and $b$ are linked)
and $a\ncap b$ means ``$a$ does not intersect $b$''
(or crossings $a$ and $b$ are not linked).
\end{defi}
 
In \cite{pos}, we gave the following 2 properties of 
Gau\ss{} diagrams, which we will use extensively in the following.

\begin{@lem}[double connectivity $2C(a,b,c)$, \cite{pos}]
\label{lem1}\rm
Whenever in a Gau\ss{} diagram $a\cap c$ and $b\cap c$
then either $a\cap b$ or there is an arrow $d$ with
$d\cap a$ and $d\cap b$.
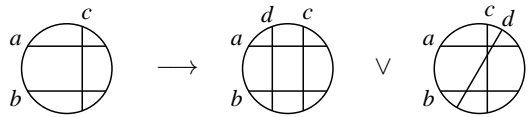
\vspace{2mm}

\[
\CD{\labch{-70}{70}{c}
\labch{30}{150}{a}
\labch{-30}{-150}{b}
} \quad\lra\quad
\CD{\labch{-70}{70}{c}
\labch{30}{150}{a}
\labch{-30}{-150}{b}
\labch{-110}{110}{d}
} \quad\lor\quad
\CD{\labch{-80}{80}{c}
\labch{30}{150}{a}
\labch{-30}{-150}{b}
\labch{240}{60}{d}
}
\]
\end{@lem}

\begin{@lem}[even valence $ev(c)$, \cite{pos}]
\label{lem2}\rm
Any chord $c$ in a Gau\ss{} diagram has odd length, i.~e., even number
of basepoints on both its sides, or equivalently,
even number of intersections with other chords.
\end{@lem}

The first step is to show the following theorem, implying that
knots with negative Fiedler invariant are not almost positive.

\begin{theorem}\label{th1}
In any almost positive diagram $K$ we have $v_3(K)\ge 0$.
\end{theorem}

\proof The idea is to show, that for any negative configuration
in the Gau\ss{} sum, that is, a configuration with a negative weight,
we can find a positive configuration, that is, a configuration with a 
positive weight, which ``equilibrates'' it. Such positive configurations
we will call accordingly ``equilibrating''.

Let $p$ be the arrow in the Gau\ss{} diagram, corresponding to
the negative crossing. There are 3 types
of negative configurations:
\def\labelenumi{\theenumi)}\mbox{}\\*[-10pt]
\begin{enumerate}
\item\label{item1} \CD{\labch{-110}{70}{b}
    \labch{-70}{110}{a}
    \labar{190}{-10}{p}
   }\quad This is equilibrated by $(a,b)$ linked.
\item\label{item2} \svCD{CD1}\addCD{\arrow{-70}{70}
		       \arrow{-110}{110}
		       \arrow{190}{-10}
		       }\drawCD{6mm}\quad 
		       By $2C(a,b,p)$, $\x c: c\cap a,c\cap b$:
\addCD{\labch{210}{30}{c}}\drawCD{6mm}\rsCD{CD1}\quad  Then $(a,b,p)$
is equilibrated by $(a,b,c)\in (4,2)0$.
\item\label{item3} \CD{\labar{195}{-15}{p}
		       \labar{165}{15}{b}
		       \labch{-90}{90}{a}
		       }\quad  This is equilibrated by $(a,b)$ linked.
\end{enumerate}
It remains to note, that no positive configuration equilibrates
this way more than one negative configuration.  \qed

To extend the result, we like to show, that, except in the
desired cases, non-equilibrating positive configurations exist,
and therefore the value of the Gau\ss{} sum is positive.
To do so, we will study the ``environment'' of the negative arrow $p$
in the Gau\ss{} diagram. In most cases we will make assumptions,
then by $2C$ and $ev$ we will show the existence of
further and further arrows in the Gau\ss{} diagram, leading at
some point unavoidably to the creation of a non-equilibrating
configuration. Therefore these assumptions turn out wrong and leave
only the desired cases. In the following we explain the tricky
details of this obvious idea.

\section{The classification}

Playing a central role in knot theory, it has been tried for a long time
to identify the unknot from its diagrams and to classify them.
Some recent progress was achieved by J.~Birman \cite{Birman}, who
developed an algorithm to recognize the unknot from its
braid representations. This algorithm, however, not unexpectedly,
is far too complex to give (or even let us hope for) some
nice explicit description of all (braid) diagrams of the unknot.
In fact, already the question which conjugacy classes of 4-braids
have unknotted closure, is known to be extremely hard \cite{%
Morton,Fied}.

More is known for special cases of diagrams. It has been proved
via different methods, that alternating \cite{Cromwell,Murasugi,%
Gabai} or positive \cite{Cromwell,pos} diagrams of the unknot
are completely reducible, that is, transformable into the
zero crossing diagram by (crossing number) reducing Reidemeister I
moves only. (One common argument is, that in such diagrams the
Seifert algorithm must give a disc, and these are exactly the
diagrams with this property.)

For almost positivity, the following appealing series of examples
comes in mind: the twist knots $3_1,4_1,5_2,6_1,7_2,8_1,\dots$,
that is, the knots with Conway notation $(k,2)$, $k\in\bN$ can be
unknotted in their alternating diagrams by 1 crossing change, giving
(modulo mirroring) an almost positive diagram, see
figure \reference{figu}(c).

Here we will show that for connected diagrams these are indeed
the only examples, which leads to a classification of all
almost positive diagrams of the unknot.

Note, that this result again gives a strong contrast to the problems
of controlling almost alternating diagrams of the unknot
\cite[\S 5.5]{Adams}.

\begin{theorem}\label{th2}
If $K$ is a connected almost positive diagram and $v_3(K)=0$,
then $K$ is an unknotted twist knot diagram or a one
crossing diagram.
\end{theorem}

\begin{figure}[htb]
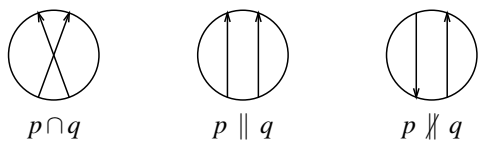

\[
\begin{array}{c@{\hspace{1cm}}c@{\hspace{1cm}}c}
  \CD{\arrow{250}{70} \arrow{290}{110}} &
  \CD{\arrow{290}{70} \arrow{250}{110}} &
  \CD{\arrow{290}{70} \arrow{110}{250}}
\\
p\cap q & p\|q & p\not\|q
\end{array}
\]
\caption{The possible mutual positions of two
arrows in a Gau\ss{} diagram.\label{fig2}}
\end{figure}

Our proof is based on some very involved analysis of the combinatorics
of the Fiedler formula, similar to this of \cite{pos}, and is divided
into several subcases. To introduce some abbreviations, 
in the following `\contr' denotes a contradiction and `$\|$'
denotes `parallel' (see figure \reference{fig2}). It appears
more appropriate to use a rather symbolic notation (even if it
reduces readability), in order to avoid misinterpretations of the
wording, as some logical constructs that will appear will be rather
complex.

\proof Fix $K$ and its negative crossing $p$. $K$ has no non-%
equilibrating  positive configurations. Therefore the following
conditions hold in $K$ (in the following we will refer to each one
by boxing its number):

\begin{enumerate}
\def\labelenumi{\theenumi.\enspace}
\item If $r,q$ are linked, $p\nin\{r,q\}$, then at least one
of $r,q$ is linked with $p$. If not both are linked with $p$,
then the not linked one is $\|p$.
\item $p$ is in \em{any} $(3,3)$ configuration. That is, whenever
$a$, $b$ and $c$ form a $(3,3)$ configuration, $p$ is one of
$a$, $b$ or $c$.
\item If $a\|b$, $p\nin\{a,b\}$
and $\x c\ne p$ with $(a,b,c)\in (4,2)0$, then $c$ is unique
with this property and $(a,b,p)\in (4,2)0$.
\item A fragment of the kind
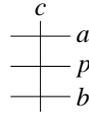
\[
\diagram{4mm}{2.5}{3.5}{
  \picline{1 0}{1 3.0}
  \picline{0 0.5}{2 0.5}
  \picline{0 1.5}{2 1.5}
  \picline{0 2.5}{2 2.5}
  \picputtext[d]{1 3.2}{$c$}
  \picputtext[l]{2.2 0.5}{$b$}
  \picputtext[l]{2.2 1.5}{$p$}
  \picputtext[l]{2.2 2.5}{$a$}
}
\]
with $c\cap a,b,p$ and $a\ncap b\ncap p\ncap a$ does not exist.
Else $a\|p$ (else $(c,a)\in\boxed1$) and $b\|p$ (else $(c,b)\in\boxed1$),
so $a\|b$ and $(a,b,c)\in (4,2)0$ \contr\ to \boxed3.
\item A fragment of the kind
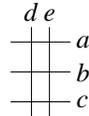
\[
\diagram{4mm}{2.5}{3.5}{
  \picline{0.7  0}{0.7 3.0}
  \picline{1.3  0}{1.3 3.0}
  \picline{0 0.5}{2 0.5}
  \picline{0 1.5}{2 1.5}
  \picline{0 2.5}{2 2.5}
  \picputtext[d]{0.7 3.2}{$d$}
  \picputtext[d]{1.3 3.2}{$e$}
  \picputtext[l]{2.2 0.5}{$c$}
  \picputtext[l]{2.2 1.5}{$b$}
  \picputtext[l]{2.2 2.5}{$a$}
}
\]
with $p\nin\{a,b,c,d,e\}$ and $e,d\cap a,b,c$ does not exist.
Else $a\cap b$ or $b\cap c$ leads to a $(3,3)$-configuration
(\contr\ to \boxed2) and for $a\ncap b$ and $b\ncap c,$
two of $\{a,b,c\}$ are $\|$, say $(a,c)$. Then $(a,c)$
participate in at least two $(4,2)0$ configurations, \contr\ to \boxed3.
\end{enumerate}

\begin{caselist}
\case $\x a,b\cap p, a\cap b$.
\[
\svCD{emptyCD}\addCD{\piclinewidth{60}\labch{190}{-10}{p}
\labch{-70}{110}{a}
\labch{-110}{70}{b}
}\drawCD{6mm}\svCD{case1}
\]
Assume,
\begin{eqn}\label{ass1}
\x c: c\ncap a,b, c\cap p\qquad
\addCD{\piccurve{40 1 x polar}{0.6 0.2}{0.6 -0.3}{-50 1 x polar}
\picputtext{40 1.2 x polar}{$c$}
}\drawCD{8mm}
\end{eqn}
Then $2C(a,c,p)\so \x d: d\cap a,d\cap c$ and 
$2C(b,c,p)\so \x d': d'\cap b,d'\cap c$. Now if $d\cap a,b$, then
$(d,a,b)\in (3,3)$ \contr, and if $d'\cap a,b$, then
$(d',a,b)\in (3,3)$ \contr, so $d'\ncap a$, $d\ncap b$ and in
particular $d\ne d'$.
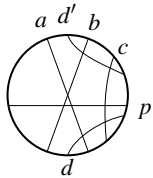
\[
\addCD{\piccurve{90 1 x polar}{0.1 0.7}{0.6 0.5}{20 1 x polar}
\picputtext[d]{90 1.2 x polar}{$d'$}
}\svCD{CD2}
\addCD{\piccurve{-20 1 x polar}{0.6 -0.4}{0.1 -0.7}{-90 1 x polar}
\picputtext{-90 1.2  x polar}{$d$}
}\drawCD{8mm}
\rsCD{CD2}
\]
Further $d\ncap d'$ (if $d\cap d'$, then by $d,d'\cap c$ we had
$(d,d',c)\in (3,3)$\contr), and therefore (see diagram) not both
$d,d'\cap p$. On the other hand, $d,d'\ncap p$ implies
$(c,d,p,d')\in\boxed 4\contr$, so
$p$ is intersected by exactly one of
$d,d'$. Assume w.l.o.g. (modulo mirroring the diagram) that
$d\cap p,d'\ncap p$.
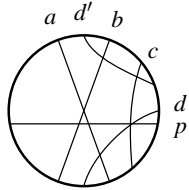
\[
\addCD{\piccurve{0 1 x polar}{0.6 -0.1}{0.1 -0.6}{-90 1 x polar}
\picputtext[dl]{0 1.2 x polar}{$d$}}
\drawCD{10mm}
\]
Because of $2C(d,c,d')$ we have $\x e\cap d,d'$. We have $e\cap p$ (else
$(e,d')\in\boxed 1$) and $e\ncap a,b,c$ (else \boxed 2). There
are 2 possibilities for $e$:
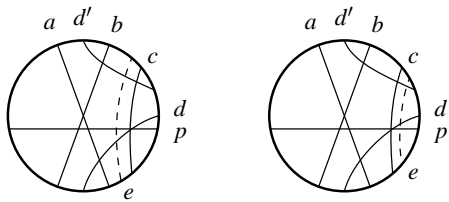
\[
\svCD{CD1}
\addCD{{\piclinedash{0.1}{0.05}
        \piccurve{50 1 x polar}{0.4 0.2}{0.4 -0.3}{-60 1 x polar}
       }
       \picputtext{-60 1.2 x polar}{$e$}
      }\drawCD{10mm}
      \rsCD{CD1}
      \hspace{10mm}
\svCD{CD1}
\addCD{{\piclinedash{0.1}{0.05}
	\piccurve{30 1 x polar}{0.75 0.2}{0.7 -0.3}{-40 1 x polar}
}
\picputtext{-40 1.2 x polar}{$e$}
}\drawCD{10mm}
\rsCD{CD1}
\]

Both choices are equivalent (the second is the same as the first
with $c,e$ swapped), so consider only the first:
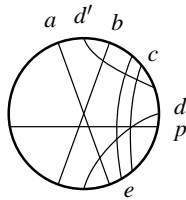
\[
\addCD{\piccurve{50 1 x polar}{0.4 0.2}{0.4 -0.3}{-60 1 x polar}
       \picputtext{-60 1.2 x polar}{$e$}
      }\drawCD{10mm}
\]
Now $2C(d',b,a)\so \x g: g\cap a,g\cap d'$. Then
$g\ncap d,e,c,b$ (else \boxed2) and $d\cap p$ (else $(g,d')\in\boxed1$).
So (modulo swapping $b$ and $g$) we obtain the following picture:
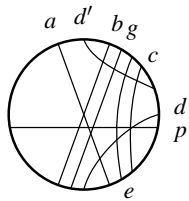
\[
\addCD{\labch{-100}{60}{g}}\drawCD{10mm}
\]
Now $ev(b)\so \x h: h\cap b$. Then $h\ncap g$ (else
$(b,g,h,a,d')\in\boxed5$) and $h\ncap a,d'$ (else \boxed2). There
are 3 choices of $h$ indicated by dashed arcs in the following picture:
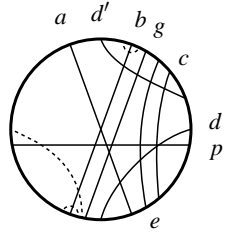
\[
\addCD{{\piclinedash{0.04}{0.02}
	\picclip{\piccircle{0 0}{1}{}}{
	  \picmultigraphics[rt]{2}{180}{
	    \piccircle{70 1 x polar}{0.09}{}
	  }
	}
	\piccurve{-102 1 x polar}{-0.2 -0.3}{-0.8 -0.1}{-1 0}
       }
      }
\drawCD{12mm}
\]
If $h$ is one of the chords in the lower part of the picture, then
apply $2C(h,b,d')$. So $\x i: i\cap h,i\cap d'$. But then $i\cap a$
and $(a,d',b,i,g)\in\boxed5$\contr. If $h$ is the chord in the upper
part of the picture, then apply $2C(h,b,a)$ to deduce $\x i:
i\cap h,i\cap a$. But then $i\cap d'$ with the same contradiction.

This shows, that the assumption \eqref{ass1} is wrong, and, using
the fact that $c\ne a,b$ with $c\cap a,b$ leads to $\boxed2$, we have
proved:
\begin{eqn}\label{eq2}
\fa c\ne a,b: c\cap p\so (c\cap a\land c\ncap b)\lor
(c\cap b\land c\ncap a)
\end{eqn}

Now assume
\begin{eqn}\label{ass2}
\x c: c\cap b,c\ncap a,c\ncap p\qquad
\rsCD{case1}\svCD{case1}
\addCD{\piccurve{90 1 x polar}{0.1 .6}{0.4 .6}{45 1 x polar}
\picputtext{45 1.3 x polar}{$c$}
}\drawCD{8mm}
\end{eqn}
The case $\x c: c\cap a,c\ncap b,c\ncap p$ is dealt with analogously.
Now $2C(a,c,b)\so \x d: d\cap c,a$. $d\ncap b$ (else \boxed2)
and $d\cap p$ (else $(d,c)\in\boxed1$). Modulo swapping $b$ and $d$
one obtains a picture like this:
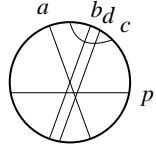
\[
\addCD{\labch{-100}{60}{d}}\drawCD{8mm}
\]
Now $ev(d)\so \x e: e\cap d$. Then $e\ncap b$ (else $(b,d,c,a,e)\in
\boxed5$) and $e\ncap c,a$ (else \boxed2). There are two choices of $e$:
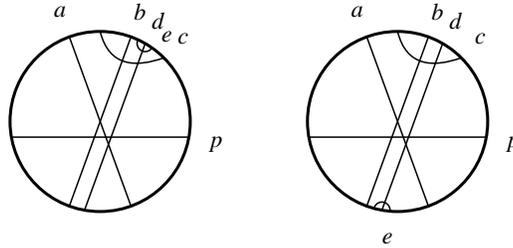
\[
\svCD{CD1}
\addCD{\picclip{\piccircle{0 0}{1}{}}{
	  \piccircle{60 1 x polar}{0.09}{}
	  }
	  \picputtext{52 1.2 x polar}{$e$}
      }\drawCD{12mm}
\rsCD{CD1}
\hspace{10mm}
\addCD{\picclip{\piccircle{0 0}{1}{}}{
	  \piccircle{-100 1 x polar}{0.09}{}
	  }
	  \picputtext{-95 1.3 x polar}{$e$}
      }\drawCD{12mm}
\]
In the first case apply $2C(e,d,a)$ to deduce $\x f: f\cap e,a$. But
then $f\cap c$ and $(a,c,b,d,f)\in\boxed5$ \contr. In the second
case apply $2C(e,d,c)$. Then $\x f: f\cap e,c$. But
then $f\cap a$ with the same contradiction.

This shows a contradiction to assumption \eqref{ass2}, so that we
obtain $(c\cap a\lor c\cap b)\so c\cap p$, and together with
\eqref{eq2} we have
\begin{eqn}
\fa c\ne a,b: (c\cap a\lor c\cap b)\iff c\cap p
\end{eqn}

Assume,
\begin{eqn}\label{ass3}
\begin{array}{c@{\ }l@{\ }l@{\quad}l@{\ }l}
& \x c:  & c\cap p, & c\cap a, & c\ncap b\\
\land & \x c': & c'\cap p, & c'\ncap a, & c'\cap b
\end{array}
\end{eqn}
\begin{caselist}
\case $c\cap c'$. Modulo swapping $(b,c)$ and $(a,c')$ we obtain
the following picture:
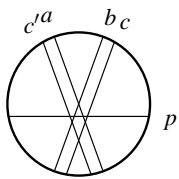
\[
\rsCD{case1}\svCD{case1}
\addCD{\labch{-100}{60}{c} \labch{-80}{120}{c'} }
\drawCD{9.5mm}
\]
Now $ev(a)\so \x d: d\cap a$. Then $d\ncap b,c$ (else a (3,3)-%
configuration of \boxed2) and $d\ncap c'$ (else
$(c',a,b,c,d)\in\boxed5$). There are three choices of $d$:
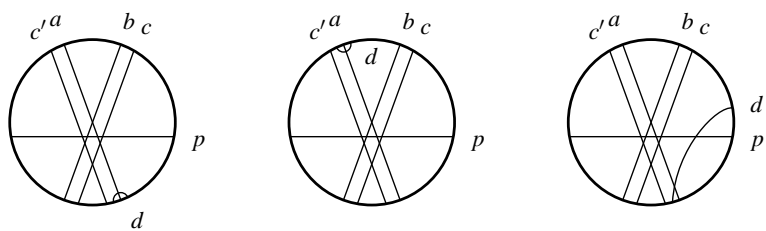
\[
\svCD{CD2}
\addCD{\picclip{\piccircle{0 0}{1}{}}{
	  \piccircle{-70 1 x polar}{0.09}{}
	  }
	  \picputtext{-65 1.3 x polar}{$d$}
      }\drawCD{11mm}
\rsCD{CD2}\hspace{10mm}
\svCD{CD2}
\addCD{\picclip{\piccircle{0 0}{1}{}}{
	  \piccircle{110 1 x polar}{0.09}{}
	  }
	  \picputtext{90 0.8 x polar}{$d$}
      }\drawCD{11mm}
\rsCD{CD2}\hspace{10mm}
\addCD{\piccurve{-75 1 x polar}{0.3 -0.3}{0.8 0.2}{10 1 x polar}
\picputtext{10 1.3 x polar}{$d$}
}\drawCD{11mm}
\]
In the first case use $2C(d,a,b)$. So $\x e: e\cap d,e\cap b$. Then
$e\cap c$ and $(b,c,a,e,c')\in\boxed5$\contr.
In the second case use $2C(d,a,c)$, so that $\x e: e\cap d,c$. Then
$e\cap b$ with the same contradiction.
In the third case use symmetry of the diagram to obtain a $d'$ like
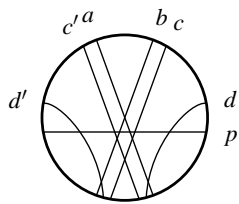
\[
\addCD{\picscale{-1 1}{
	  \piccurve{-75 1 x polar}{0.3 -0.3}{0.8 0.2}{10 1 x polar}
       }
       \picputtext{170 1.3 x polar}{$d'$}
}\drawCD{11mm}
\]
(or applying one of the other two cases to $d'$ to obtain
a contradiction). Then $2C(d,d',p)\so \x e:e\cap d,e\cap d'$.
But then $e\cap c,e\cap c'$, which together with $c\cap c'$ implies
$(c,c',e)\in(3,3)$\contr\ to \boxed2.

\case $c\ncap c'$.
\[
\rsCD{case1}
\addCD{\picmultigraphics[S]{2}{-1 1}
	  {\piccurve{-82 1 x polar}{0.1 -0.1}{0.2 0.3}{60 1 x polar}}
       \picputtext{60 1.3 x polar}{$c$}
       \picputtext{120 1.3 x polar}{$c'$}
      }\drawCD{10mm}
\]
Because of $2C(c',p,c)$ we have $\x d: d\cap c,c'$. Then $d\ncap a,b$
(else \boxed2). So the diagram looks like:
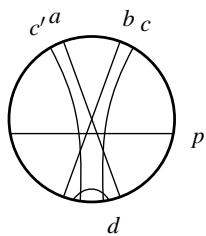
\[
\addCD{\piccurve{-78 1 x polar}{0.1 -0.8}{-0.1 -0.8}{-103 1 x polar}
\picputtext{-78 1.3 x polar}{$d$}
}\drawCD{11mm}
\]
Now $2C(d,c',b)\so \x e: e\cap d,e\cap b$. Then $e\cap p$ (else
$(e,d)\in\boxed1$) and $e\ncap c,c',a$ (else \boxed2). In the same way,
$2C(d,c,a)\so \x f: f\cap a,d$, and $f\cap p$, $f\ncap c,c',b$.
Furthermore, $f\ncap e$ (else $(f,e,d)\in(3,3)$). Modulo swapping
$(f,c)$ and $(e,c')$ we obtain the following picture:
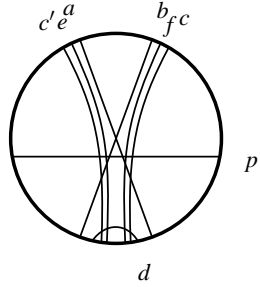
\[
\addCD{\picmultigraphics[S]{2}{-1 1}
	   {\piccurve{-85 1 x polar}{0.05 -0.1}{0.1 0.2}{65 1 x polar}}
       \picputtext{65 1.2 x polar}{$f$}
       \picputtext{115 1.2 x polar}{$e$}
}\drawCD{14mm}
\]
Now $ev(e)\so \x g: g\cap e$. For this $g$, $g\ncap d,b$ (else
\boxed2), $g\ncap c'$ (else $(c',e,b,d,g)\in\boxed5$) and $g\cap p$
(else $(g,p,d,e)\in\boxed4$). However, a glimpse at the above diagram
shows, that no such $g$ exists. Therefore, this contradiction
shows, that our assumption \eqref{ass3} is wrong, and we have
\begin{eqn}\label{eq4}
\x i\in\{a,b\}: \fa c\ne a,b: p\cap c\iff p\cap i\land p\ncap\bar i
\end{eqn}
with $\bar i\in\{a,b\}\setminus \{i\}$. Assume w.l.o.g. $i=a$.
As a consequence of \eqref{eq4}, all $c\ne a: p\cap c$ do not mutually
intersect (else an intersecting pair would build a $(3,3)$-%
configuration with $a$). 
Then the subdiagram made up of $p$ and all its neighbors, i.e.,
all crossings linked with $p$, looks like:
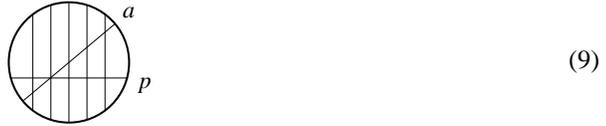
\begin{eqn}\label{GD}
\rsCD{emptyCD}\svCD{emptyCD}
\addCD{\piclinewidth{60}\picclip{\piccircle{0 0}{1}{}}
	       {\picmultigraphics{5}{0.3 0}{\picline{-0.6 -1}{-0.6 1}}}
	       \labch{195}{-15}{p}
	       \labch{-140}{40}{a}
}\drawCD{8mm}
\end{eqn}
Assume now, there are more chords in the Gau\ss{} diagram. So by
connectedness
$\x d:d\cap c$ for some $c\cap p$ and $d\ncap p$. So in particular
$c\ne a$ (as $d\cap a\so d\cap p$). 
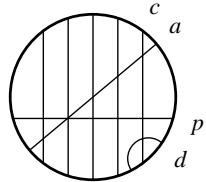
\[
\addCD{\picclip{\piccircle{0 0}{1}{}}
	  {\piccircle{-48 1 x polar}{0.25}{}}
	   \picputtext{-35 1.3 x polar}{$d$}
	   \picputtext{55 1.3 x polar}{$c$}
}\drawCD{11mm}
\]
But then $ev(c)\so \x d'\cap c$. $d'\ncap d$ (else \boxed2) and
by assumption $d,d'\ncap p$, so $(c,d,p,d')\in\boxed4$\contr.
Therefore there are no more chords in the Gau\ss{} diagram than those of
\eqref{GD}, and the knot diagram is an unknotted odd crossing number
twist knot diagram.
\end{caselist}

\case No pair $a,b$ with $a,b\cap p$ is linked. Then the subdiagram 
of the Gau\ss{} diagram made up of $p$ and all crossings linked with
it looks like:
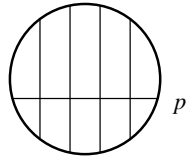
\[
\rsCD{emptyCD}\svCD{emptyCD}
\addCD{\piclinewidth{60}\picclip{\piccircle{0 0}{1}{}}
	       {\picmultigraphics{4}{0.4 0}{\picline{-0.6 -1}{-0.6 1}}}
}\svCD{even2CD}
\addCD{\labch{195}{-15}{p}}
\drawCD{1cm}
\]
By the assumption and exclusion of \boxed2 and \boxed4, no
chord intersecting $p$ is intersected more than once else.
By even valence, then any $c: c\cap p$ must be intersected exactly
once else.

Assume now $a,b\cap p$. Then by $2C(a,b,p)$ the second
chord intersecting $a$ and $b$ must be the same for all $a,b\cap p$.

But then the Gau\ss{} diagram looks like:
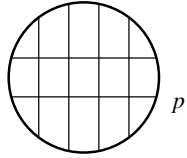
\[
\addCD{\chrd{165}{15}}\drawCD{1cm}
\svCD{CDeven}
\]
and the knot diagram is an unknotted even crossing number
twist knot diagram. \qed
\end{caselist}

An immediate consequence of the theorem is

\begin{corr}\label{Cr3.1}
No knots with $v_3=0$ (\em{inter alia}, achiral knots) can be
almost positive. \qed
\end{corr}

This result seems to have been first obtained (without published proof)
by Przytycki and Taniyama \cite{PrzTan}, and has been recently recovered
in a nice way by Lee Rudolph \cite{Rudolph}.

\begin{corr}\label{Cr3.2}
All almost positive unknot diagrams are connected sums of diagrams as in
figure \reference{figu} with summands of the kind (b) and (c) appearing
together exactly once. \qed
\end{corr}

\begin{figure}[htb]
\[
\begin{array}{c@{\hspace{1cm}}c@{\hspace{1cm}}c}
\diag{8mm}{1}{3}{
  \picPSgraphics{0 setlinecap}
  \rbraid{0.5 1.5}{1 2}
  \piccirclearc{0.5 0.5}{0.5}{-180 0}
  \piccirclearc{0.5 2.5}{0.5}{0 180}
  }
   &
\diag{8mm}{1}{3}{
  \picPSgraphics{0 setlinecap}
  \lbraid{0.5 1.5}{1 2}
  \piccirclearc{0.5 0.5}{0.5}{-180 0}
  \piccirclearc{0.5 2.5}{0.5}{0 180}
} &
\diag{5mm}{12}{6}{
  \pictranslate{0 0.54}{
  \picPSgraphics{0 setlinecap}
  \picrotate{-90}{
     \rbraid{-0.5 3}{1 2}
     \rbraid{-0.5 5}{1 2}
     \rbraid{-0.5 9}{1 2}
     \picmultigraphics{3}{0 0.4}{\picputtext[d]{-0.5 6.6}{.}}
  }
     \pictranslate{6 0}{
       \picmultigraphics[S]{2}{-1 1}{
	 \piccirclearc{4 2.5}{1.5}{-90 90}
	 \piccirclearc{4 2.5}{2.5}{-90 90}
	 \picline{4 4}{0.2 4}
	 \picline{4 5}{0.2 5}
	 \picmulticirclearc{-4.9 1 -1.0 0}{0.2 4.5}{0.5}{90 270}
       }
     }
 }
}
\\[3mm]
(a) & (b) & (c)
\end{array}
\]
\caption{\label{figu}}
\end{figure}

\section{Whitehead doubles and the Casson invariant\label{SCas}}


The use of Vassiliev invariants allows to extend the previous chirality
results also to certain cables of positive and almost positive knots.

\begin{defi} 
In the following picture we summerize three kinds of clasps in
a knot diagram and how we will call them (the diagrams are understood
up to mirroring and strand orientation, when latter is not specified).
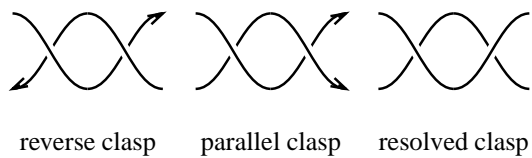
\[
\begin{tabular}{c@{\quad}c@{\quad}c}
\diag{1cm}{2}{1}{
  \picrotate{-90}{
    \lbraid{-0.5 0.5}{1 1}
    \lbraid{-0.5 1.5}{1 1}
    \picvecline{-0.95 1.9}{-1 2}
    \picvecline{-0.05 0.1}{0 0}
  }
} &
\diag{1cm}{2}{1}{
  \picrotate{-90}{
    \lbraid{-0.5 0.5}{1 1}
    \lbraid{-0.5 1.5}{1 1}
    \picvecline{-0.95 1.9}{-1 2}
    \picvecline{-0.05 1.9}{0 2}
  }
} &
\diag{1cm}{2}{1}{
  \picrotate{-90}{
    \lbraid{-0.5 0.5}{1 1}
    \rbraid{-0.5 1.5}{1 1}
  }
} \\[8mm]
reverse clasp & parallel clasp & resolved clasp
\end{tabular}
\]
Formally, a clasp is a digon, a connected component of the complement
of the diagram, neighboring just two crossings. It can be
identified with the (unordered) pair of these two crossings.
\end{defi} 

The idea to consider untwisted Whitehead doubles came out of a nice
relationship between the Vassiliev invariants of degrees 2 and 3. Let
$v_2$ denote the Vassiliev invariant of degree 2 given by
\[
v_2\,=\,-\frac 16 V''(1)\,=\,[\nabla(z)]_{z^2}
\]
(in the following $\Delta$ is the Alexander polynomial,
$\nabla$ the Conway polynomial \cite{Conway} and $P$ the HOMFLY
polynomial)
and $w_\pm$ denote the untwisted double with positive (resp. negative)
clasp. Then in \cite{pos} it was proved that
\begin{eqn}\label{v23w}
v_3\left(w_{\pm}(K)\right)=\pm 8 v_2(K)\,.
\end{eqn}
This follows basically from a combination of two observations of
Rong-McDaniel \cite{McDanielRong} (that the dualization of a
double operation of a knot maps Vassiliev invariants to Vassiliev
invariants) and Lin \cite{Lin} (that this endomorphism is nilpotent).

Whitehead doubles are classical examples of knots with $\Delta=1$
and so have been under some consideration in connection with the often
raised question about a non-trivial knot $K$ with $V_K=1$ (or even
$P_K=1$). In \cite{LickThis} Lickorish and Thistlethwaite excluded
Whitehead doubles of adequate knots from having this property, and in
\cite{Rudolph2}, as quoted in \cite{KalLin}, it was shown that for
$K$ positive $P(w_{\pm}(K))\ne 1$. Latter result inspired me to
give an strengthening of it in \cite{pos}
by combining \eqref{v23w} with the lower bound for $v_2$ in positive
diagrams in \cite{pos}, showing
that in fact for $K$ positive $w_{\pm}(K)$ has non-zero
degree 3 Vassiliev invariant, and hence in particular is chiral
and has non-trivial Jones polynomial. The above corollary then
follows from a similar inequality to this of $v_3$ in theorem
\reference{th2}.

Note, that for Whitehead doubles neither the signature, nor
the Bennequin inequality work (in some easy way) to show chirality.

Unfortunately, the need to exclude the knotted cases forces us
to reprove, this time using $v_2$, the result on almost positive
unknot diagrams. Fortunately, this time the proof is somewhat easier.

\begin{theorem}\label{th3}
If $D$ is a connected almost positive diagram then $v_2(D)\ge 0$,
and if $v_2(D)=0$, then $D$ is an unknotted twist knot diagram or a one
crossing diagram.
\end{theorem}

An immediate corollary of this theorem generalizes a recent result of
Menasco and Zhang \cite{MenZh} on the property $P$ conjecture
\cite{BM} (see e.g. \cite{BS,DR,CochranGompf}).

\begin{corr}
Almost positive knots have property $P$ (i.e. every nontrivial surgery
on $S^3$ along such knots produces a non-simply connected manifold).
\end{corr}

\proof This follows from the surgery formula for the Casson invariant
of homology spheres \cite{Akbulut} (see also \cite[\S 1]{MenZh}). \qed

For the proof of theorem \reference{th3} we need some preparations.
Beside $ev$ and $2C$ we will need the following elementary observation.

\begin{@lem}[extended even valence $eev(c)$, see \cite{pos}]\rm
In the Gau\ss{} diagram of a positive knot diagram,
exactly one half of the arrows intersecting any chord $c$
intersect it in one or the other direction (that is, are distinguished
or not in the resulting linked pair). \qed
\end{@lem}

We will need the following Gau\ss{} sum formula for $v_2$ which is due
to Polyak and Viro \cite[(3)]{VirPol}:
\begin{eqn}\label{q1}
v_2\,=\,\GD{\arrow{225}{45}
\arrow{315}{135}
\picfillgraycol{0}\picfilledcircle{1 90 polar}{0.09}{}
}\,=\,
\frac 12\left(\GD{\arrow{225}{45}
\arrow{315}{135}
\picfillgraycol{0}\picfilledcircle{1 90 polar}{0.09}{}
} + \GD{\arrow{45}{225}
\arrow{135}{315}
\picfillgraycol{0}\picfilledcircle{1 90 polar}{0.09}{}
} \right)\,.
\end{eqn}
The point on the circle depicts a point to be put somewhere
on the knot curve in the diagram, but not at a crossing.

\proof[of theorem \reference{th3}]
The fact that $v_2(D)\ge 0$ is almost positive diagrams may
be provable along similar lines as $v_3$, but it is easier to quote
Cromwell's (skein-theoretic) proof \cite[corollary 2.2, p.~539]
{Cromwell}.

Now consider an almost positive knot diagram $D$ and its Gau\ss{}
diagram with the negative arrow $k$. Recall the move we used in the
proof of the positivity of $v_2$ in positive diagrams in \cite{pos}:
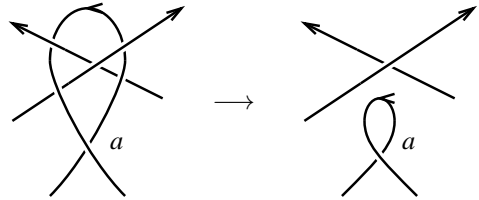
\[
\diag{1cm}{3}{3.5}{
  \pictranslate{-1 0}{
    \picPSgraphics{0 setlinejoin}
    \picline{1.5 1.5}{2.25 2}
    \picline{3.5 1.8}{2.2 2.45}
    \picmulticurve{0.12 1 -1.0 0}{2 0.5}{2.5 1}{3 2}{3 2.3}
    \picmulticurve{0.12 1 -1.0 0}{3 0.5}{2.5 1}{2 2}{2 2.3}
    \picellipsevecarc{2.5 2.3}{0.5 0.7}{0 90}
    \picellipsearc{2.5 2.3}{0.5 0.7}{90 180}
    \picmultivecline{0.12 1 -1.0 0}{2.2 2.45}{1.5 2.8}
    \picmultivecline{0.12 1 -1.0 0}{2.25 2}{3.75 3}
    \picputtext[l]{2.8 1.2}{$a$}
  }
}\,\lra\,
\diag{1cm}{3}{3.5}{
  \pictranslate{-1 0}{
    \picPSgraphics{0 setlinejoin}
    \picline{1.5 1.5}{2.25 2}
    \picline{3.5 1.8}{2.2 2.45}
    \picmulticurve{0.12 1 -1.0 0}{2 0.5}{2.5 1}{2.7 1.2}{2.7 1.5}
    \picmulticurve{0.12 1 -1.0 0}{3 0.5}{2.5 1}{2.3 1.2}{2.3 1.5}
    \picellipsevecarc{2.5 1.5}{0.2 0.3}{0 90}
    \picellipsearc{2.5 1.5}{0.2 0.3}{90 180}
    \picmultivecline{0.12 1 -1.0 0}{2.2 2.45}{1.5 2.8}
    \picmultivecline{0.12 1 -1.0 0}{2.25 2}{3.75 3}
    \picputtext[l]{2.8 1.2}{$a$}
  }
}\,.
\]
We trivialized loops by pulling them above the rest of the diagram by
crossing changes. On the level of Gau\ss{} diagrams this means
that we delete a chord $a$, intersecting all chords ending on one of its
sides, and then also delete all these intersecting chords.

For the following arguments it is convenient to place the point
of  \eqref{q1} near one of the endpoints of $a$. Now consider what
happens with $v_2$ under our move.

It is clear that if $a\ncap k$, then the move never augments $v_2$,
because the only negative (contribution) configurations removed are
those of $(k,c)$ with $c\cap a$, but their contribution is equilibrated
by those of $(c,a)$. Moreover, because of $ev$ and $eev$ an even number
of $c$ with $c\cap a$ intersects $k$ (because after the move the
Gau\ss{} diagram still corresponds to a knot diagram), an even
number of them does not intersect $k$, and exactly half of these numbers
intersect $c$ in either direction.
Hence if $a\cap c$ with $c\ncap k$, resolving $a$ would strictly reduce
$v_2$ (as $a$ is linked with some $c$ by connectedness) 
and leave over an almost positive diagram. Therefore, after the
move still $v_2\ge 0$ and hence before the move $v_2>0$. So for any
$a$, $a\cap k$ or $a$ intersects all $c\cap k$. This splits the
positive arrows in two parts
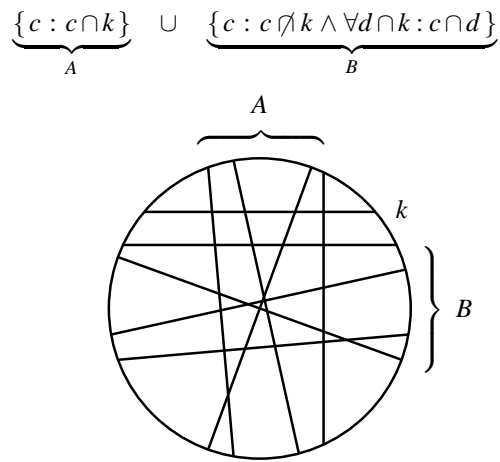
\[
\underbrace{\{\,c\,:\,c\cap k\,\}}_{A}\quad\cup\quad
\underbrace{\{\,c\,:\,c\ncap k\,\land\,\fa d\cap k: c\cap d\,\}}_{B}
\]
\[
\diag{1cm}{5}{5}{
\pictranslate{2.56 2.56}{
  \picmultigraphics[rt]{4}{90}{
    \piccirclearc{0 0}{2}{-45 45}
  }
  \picputtext[d]{-0.0 2.6}{$A$}
  \picputtext[l]{2.6 0 polar}{$B$}
  \picmultigraphics[rt]{2}{90}{
    \picline{-65 2 x polar}{65 2 x polar}
    \picline{-110 2 x polar}{70 2 x polar}
    \picline{-75 2 x polar}{100 2 x polar}
    \picline{110 2 x polar}{-100 2 x polar}
    \picputtext{2.3 10 polar}{$\left.\rule{0mm}{1cm}\right\}$}
  }
  \picline{2 40 polar}{2 140 polar}
  \picputtext{2.3 35 polar}{$k$}
}
}
\]
Now putting the point on some segment of the circle between an
endpoint of a chord in $A$ and an endpoint of a chord in $B$
and using $ev$ and $eev$ we see that any $c\in B$ gives a positive
contribution to $v_2$ equilibrating the negative contribution of $k$.
Hence for $v_2=0$ we must have $|B|\le 1$.

If $B=\{c\}$, then $c$ forms with $k$ an (unlinked) resolved
clasp, whose elimination by a Reidemeister II move gives a positive
diagram with $v_2=0$. This diagram must then have only isolated
chords. Hence $A$ has only non-mutually intersecting arrows
(or non-linked crossings in the knot diagram), whose number by $ev$ must
be even, so we have
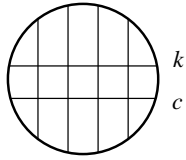
\[
\rsCD{even2CD}\addCD{\labch{195}{-15}{c}
\chrd{10}{170}\picputtext{1.3 12 polar}{$k$}
}
\drawCD{1cm}
\]
an unknotted twist knot diagram of even crossing number.

If $B=\varnothing$, then the knot diagram can be obtained just by
Reidemeister II moves and crossing changes from
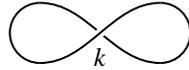
\[
\diag{8mm}{3}{1}{
\pictranslate{1.5 0.5}{
  \picrotate{90}{
    \picPSgraphics{0 setlinecap}
    \pictranslate{-0.5 -1.5}{
    \rbraid{0.5 1.5}{1 2}
    \piccirclearc{0.5 0.5}{0.5}{-180 0}
    \piccirclearc{0.5 2.5}{0.5}{0 180}
    }
  }
  \picputtext{0 -0.4}{$k$}
}}
\]
But then clearly the diagram will have a clasp somewhere
(at least the one created by the last Reidemeister II move). Hence
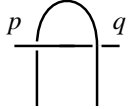
\[
\diag{1cm}{2}{1.6}{
  \picline{1 1}{1.7 1}
  \picmulticurve{0.12 1 -1.0 0}{1.4 1}{1.4 1.8}{0.6 1.8}{0.6 1}
  \picmultiline{0.12 1 -1.0 0}{1.4 1}{1.4 0.2}
  \picline{0.6 0.2}{0.6 1}
  \picmultiline{0.12 1 -1.0 0}{0.3 1}{1 1}
  \picline{0.9 1}{1.1 1}
  \picputtext{1.7 1.3}{$q$}
  \picputtext{0.3 1.3}{$p$}
}
\]
all positive crossings in the diagram can be resolved by
consecutively removing clasps. Putting the point near one of the
endpoints of $k$ in the Gau\ss{} diagram, one sees because of $ev$
and $eev$ that resolving such a clasp $(p,q)$ never augments $v_2$
and that it reduces it strictly if $\x c\ne k$ with $c\cap p,q$
and $p\cap q$ or if $\x c_1,c_2\ne k$ with $c_1,c_2\cap p,q$ and
$p\ncap q$.
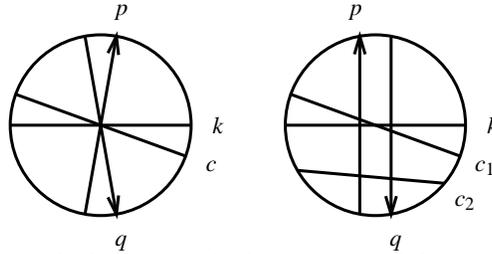
\[
\szCD{1.2cm}{
  \labar{-100}{80}{p}
  \labar{100}{-80}{q}
  \labch{-180}{0}{k}
  \labch{-200}{-20}{c}
} \qquad
\szCD{1.2cm}{
  \labar{-100}{100}{p}
  \labar{80}{-80}{q}
  \labch{-180}{0}{k}
  \labch{-200}{-20}{c_1}
  \labch{-150}{-40}{c_2}
}
\]
This means that for $v_2=0$ no such clasp occurs in
the sequence of clasps to be resolved. In particular there cannot be
$a,b,c\cap k$ with $a\cap b,c$ and $b\cap c$,
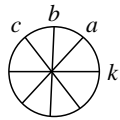
\[
\GD{
  \labch{-180}{0}{k}
  \labch{-135}{50}{a}
  \labch{-95}{90}{b}
  \labch{-55}{130}{c}
}
\]
as at some point $a$, $b$ or $c$ must be involved into a clasp
and it would have to be one of the above kinds. Consequently,
no subdiagram like
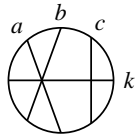
\begin{eqn}\label{**}
\rsCD{emptyCD}\addCD{
  \labch{-180}{0}{k}
  \labch{-55}{55}{c}
  \labch{-130}{90}{b}
  \labch{-90}{130}{a}
}\drawCD{7mm}\svCD{CD2}
\end{eqn}
occurs in the Gau\ss{} diagram.
Else by $2C(b,k,c)$, $\x d\cap b,c$. By assumption $d\cap k$, so
$d\ncap a$ (else $d,a,b\in\eqref{**}$). But then
by $2C(a,k,c)$ similarly there is some $d'$
such that $d,d'$ and $c$ mutually intersect (contradicting \eqref{**}).
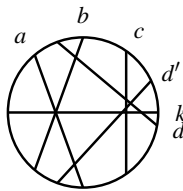
\[
\rsCD{CD2}
\addCD{\labch{110}{-10}{d}\labch{-110}{25}{d'}}
\drawCD{1cm}
\]
So if $a, b\cap k$ and $a\cap b$, any $c\cap k$ intersects
exactly one of $a$ and $b$. This pairs up the chords intersecting $k$
in 2 parts of $p$ and $q$ elements, such that $a\cap b$ iff $a$ and
$b$ belong to different parts. $ev$ forces $p$ and $q$ to be odd,
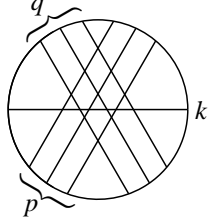
\[
\diag{6mm}{5}{5}{
\pictranslate{2.56 2.56}{
  \picmultigraphics[rt]{2}{120}{
    \piccirclearc{0 0}{2}{0 240}
    \picrotate{30}{
    \picline{-70 2 x polar}{70 2 x polar}
    \picline{-85 2 x polar}{85 2 x polar}
    \picline{-100 2 x polar}{100 2 x polar}
    \picrotate{90}{
      \picputtext{2.2 5 polar}{$\left.\rule{0mm}{0.5cm}\right\}$}
    }
    }
  }
  \picputtext{2.7 120 polar}{$q$}
  \picputtext{2.6 -125 polar}{$p$}
  \picline{2 180 polar}{2 0 polar}
  \picputtext{2.3 0 polar}{$k$}
}
}
\]
and then the diagram is the diagram of a $(p,q,-1)$ pretzel knot
with $p,q>0$ odd. Finally, to show $v_2=0$ only if such a knot
is unknotted is a matter of direct calculation. It is known (or
can be deduced from the formula for $v_2$) that $v_2(P(p,q,r))=
\mbox{\small$\ds\frac{pq+pr+qr+1}{4}$}$, so $v_2(P(p,q,-1))=0$
implies $p=1$ or $q=1$, and we have an unknotted twist knot
diagram of odd crossing number. \qed

\begin{corr}\label{Cr4.1}
The  untwisted Whitehead doubles (with either clasps) of an
almost positive knot are chiral, and have non-trivial Jones
polynomial. 
\end{corr}

\proof
As noted, this is straightforward
from theorem \reference{th3} and \eqref{v23w}. \qed

Clearly the chirality of the satellite follows from that of the
companion, hence this result is a consequence of the previous one
(corollary \ref{Cr3.1})
also by classical arguments. However, the following corollary
shows that the chirality result extends to many cases, where it is
less obvious, for example for the connected sum of an almost positive
knot and its obverse.

\begin{corr}
If a knot $K$ is the connected sum of positive, almost positive
knots and their obverses, or a cable knot thereof, then
its untwisted Whitehead doubles (with either clasps) 
are chiral (and have non-trivial Jones polynomial).
\end{corr}

\proof If $K$ is a connected sum, then the positivity of $v_2$
follows from its invariance under mirroring and additivity under
connected sum. As for cables, if $p,q\in\bN$ are coprime, $T_{p,q}$
denotes the $(p,q)$-torus knot and $K_{p,q}$ the satellite of $T_{p,q}$
around $K$, we have
\begin{eqn}\label{sat}
v_2(K_{p,q})\,=\,a_{p,q}\,v_2(K)+c_{p,q}\,.
\end{eqn}
It is known that the Eigenvalues of Vassiliev invariants (modulo
Vassiliev invariants of lower degree) under cabling operations
are always positive (see \cite{KSA} or \cite{McDanielRong}; in fact,
they are powers of the number of parallels of the cable), hence
so is $a_{p,q}$, and putting the unknot in \eqref{sat},
we obtain $c_{p,q}=v_2(T_{p,q})>0$, hence the positivity of $v_2$
for the cable follows from that of its companion. \qed

Finally, we can now exhibit non-triviality of the polynomials
for almost positive knots themselves. While for $V$ (and hence for $P$
and $F$) we could have concluded non-triviality already in the previous
paragraph using $v_3$, we postponed it till now in order
to do it for all five polynomials at one go. This is another
considerable bonus of using Vassiliev invariants instead of signatures
as in \cite{PrzTan}.

\begin{corr}
For any of the polynomials of Alexander/Conway, Jones, HOMFLY,
Brandt-Lickorish-Millett-Ho and Kauffman there is no almost
positive knot with unit polynomial.
\end{corr}

\proof Use the positivity of $v_2$ and the relations
\[
-6v_2\,:=-3\Delta''(1)=V''(1)=Q'(-2)
\]
and the well-known specializations for the HOMFLY and Kauffman
polynomial. The equality between the Jones and Alexander polynomial is
probably due already to Jones \cite[\S 12]{Jones2}. The relation
between the Jones and Brandt-Lickorish-Millett-Ho polynomial is
proved by Kanenobu in \cite{Kanenobu}. \qed

For the Jones (and hence also HOMFLY and Kauffman) polynomial,
this result was announced also in \cite{PrzTan}. An alternative
proof, using arguments similar to those of Przytycki and Taniyama,
appeared in \cite{restr}.

\section{Some further inequalities related to the genus\label{Sgen}}

An obvious desire is to improve the positivity results for $v_2$
and $v_3$ on almost positive knots to inequalities involving the
crossing number, as for positive knots in \cite[theorem 6.1]{pos}.
Focusing in the following on $v_2$, which leads to more interesting
consequences, this is related to a conjecture made in \cite[\S 6]{pos}.

\begin{defi}
Call a diagram \em{bireduced}, if it is reduced, i.e. has no nugatory
crossings, and does not admit a move
\begin{eqn}\label{second}
\diag{1cm}{2}{1.6}{
  \picline{0.3 1}{1 1}
  \picmulticurve{0.12 1 -1.0 0}{1.0 1.5}{0.7 1.5}{0.4 1.1}{0.7 0.8}
  \picmulticurve{0.12 1 -1.0 0}{1.3 0.8}{1.6 1.1}{1.3 1.5}{1.0 1.5}
  \picmultiline{0.12 1 -1.0 0}{1 1}{1.7 1}
  \picline{0.7 0.2}{1.3 0.8}
  \picmultiline{0.12 1 -1.0 0}{0.7 0.8}{1.3 0.2}
  \piccurve{1.0 1.5}{0.7 1.5}{0.4 1.1}{0.7 0.8}
}\quad\lra\quad
\diag{1cm}{2}{1.6}{
  \picline{1 1}{1.7 1}
  \picmulticurve{0.12 1 -1.0 0}{1.4 1}{1.4 1.8}{0.6 1.8}{0.6 1}
  \picmultiline{0.12 1 -1.0 0}{1.4 1}{1.4 0.2}
  \picline{0.6 0.2}{0.6 1}
  \picmultiline{0.12 1 -1.0 0}{0.3 1}{1 1}
  \picline{0.9 1}{1.1 1}
}\,.
\end{eqn}
\end{defi}

{}From now on we will assume that all diagrams are reduced.

\begin{conj}\label{cj3}
In a positive bireduced diagram $D$, $v_2\ge lk(D)/4$, where $lk(D)$
is the number of linked pairs in $D$.
\end{conj}

A computer experiment revealed that the conjecture is in general
false. One counterexample is 
a diagram of $8_{19}$ obtained by making positive by crossing changes
the alternating diagram of $9_{40}$. It has 21 linked pairs,
but $v_2=v_2(8_{19})=5$.

Nevertheless, a weaker version of the conjecture is true.
To formulate the statement, we need to recall the notion of the
genus of a diagram. In the sequel we will also use the weak genus
of a knot $K$. Both terms were introduced in \cite{gen1}.

\begin{defi}
For a diagram $D$ of knot $K$, we define the genus $g(D)$ as the
genus of the surface obtained by applying the Seifert algorithm to
this diagram:
\[
g(D)=\frac{c(D)-s(D)+1}{2}\,,
\]
with $c(D)$ and $s(D)$ being the crossing and Seifert circle number of
$D$, respectively. The weak genus of $K$, denoted by $\tl g(K)$, is
the minimal genus of all its diagrams:
\[
\tl g(K)\,=\,\,\min\left\{\,g(D)\,:\,\mbox{$D$ is
a diagram of $K$}\,\right\}\,.
\]
\end{defi}

\begin{theo}\label{T3}
For every genus $g$, and every $\eps>0$, there are at most finitely
many positive diagrams $D$ of genus $g$ and $lk(D)/v_2(D)>4+\eps$.
\end{theo}

We recall the main result of our work on the diagram genus in
\cite{gen1}, which will be used in the following.
It is related to the move below which we call a $\bt$ move:
\begin{eqn}\label{t2m}
\diag{7mm}{1}{1}{
    \picmultivecline{0.12 1 -1.0 0}{0 0}{1 1}
    \picmultivecline{0.12 1 -1.0 0}{1 0}{0 1}
}\lra\quad
\diag{7mm}{3}{2}{
  \picPSgraphics{0 setlinecap}
  \pictranslate{0.5 1}{
    \picrotate{-90}{
      \lbraid{0 -0.5}{1 1}
      \lbraid{0 0.5}{1 1}
      \lbraid{0 1.5}{1 1}
      \pictranslate{-0.5 0}{
      \picvecline{0.03 1.95}{0 2}
      \picvecline{0.03 -.95}{0 -1}
    }
    }
    }
}
\,.
\end{eqn}
Another move we will use, the \em{flype}, is the one of \cite{MenThis}.

\begin{theo}(theorem  3.1 of \cite{gen1})\label{T4} Reduced knot
diagrams of given genus decompose into finitely many equivalence classes
modulo crossing changes, $\bt$ moves and their inverses. That
is, reduced knot diagrams of given genus with no resolved clasps
can be all obtained from finitely many (called ``generating'') diagrams
by repeated $\bt$ moves.
\end{theo}

A further result we will need later is proved in \cite{pos}.

\begin{lem}\label{LM2}(\cite{pos})
If $D$ is a positive bireduced diagram of $c$ crossings, then
$lk(D)\ge 3\mybrtwo{c-1}$.
\end{lem}

\proof[of theorem \ref{T3}]
One needs to observe that if $D'$ arises from $D$ by a $\bt$ move, then
$lk(D')-lk(D)=4\bigl(v_2(D')-v_2(D)\bigr)$. To see this, put
the basepoint in \eqref{q1} near the (crossings of the) created clasp.
Then theorem \ref{T4} implies that for
any sequence $D_1,D_2,\dots$ of (distinct) positive diagrams of fixed
genus, it holds $lk(D_i)/v_2(D_i) \to 4$ as $i\to\infty$. \qed

The number of generating diagrams grows very rapidly with
the genus, but in \cite{gen1,gen2} at least a description of all such
diagrams for genus one and two was obtained. This work can be used to
make the estimate in theorem \reference{T3} more explicit
for small genus by accounting for the exceptional cases.

\begin{prop}\label{PP1}
On any positive diagram of genus at most $3$, we have $v_2/lk\ge 5/21$.
\end{prop}

\proof
To make a systematic verification of the maximal ratio $lk/v_2$ on
positive diagrams, once having it found to be $>4$, the argument
proving theorem \ref{T3} shows that we need to consider just
diagrams without reverse clasps. A similar argument shows the same also
for parallel clasps (resolving a parallel clasp in a positive diagram
reduces $lk$ by $4n-1$ and $v_2$ by $n$). Thus let $D$ be a positive
diagram of genus at most 3 with no clasps.

We fix a linked pair $(a,b)$ of crossings in $D$ and smooth them out,
obtaining a genus 2 diagram $D'$. If $D$ cannot be simplified (after
a possible sequence of flypes) by the inverse of a $\bt$ move,
we showed in the proof of theorem 3.1 of \cite{gen1}, that $D'$
has at most 4 reducible crossings, but using the stronger condition
that $D$ has no clasp, we see that out argument there
modifies to show that in fact there is at most one 
reducible crossing $p$ in $D'$. Let $D''$ be $D'$ if $p$ does not exist,
or the diagram obtained from $D'$ after reducing $p$ according
to \eqref{eqred}, if $p$ exists. 

It is easy to see that for \em{knot} diagrams with $p$
being the only reducible crossing on the left of \eqref{eqred},
reducing $p$ never augments the number of clasps by more than
one\footnote{Note, that this is not true for link diagrams~-- consider
the diagram of the Hopf link with one kink.}. Also, as smoothing out
a crossing augments the number of clasps at most by 2, $D'$ has
at most $4$ clasps (note that one crossing may be in two clasps,
as on the right of \eqref{t2m}).

\begin{caselist}
\case
If $D''$ is connected and $p$ existed, then one of $P$ and $Q$
in \eqref{eqred} has no crossing. In this case, $a$, $b$ and $p$
were the corners of a triangular component of the complement of $D$,
and smoothing out $b$ destroyed (one of) the clasp(s) created by
smoothing out $a$. Thus $D'$ has (instead of at most $4$) in fact
at most $2$ clasps, and $D''$ has at most $3$ clasps.

\case
If $D''$ is connected and $p$ did not exist, then $D''=D'$ has
at most $4$ clasps, as noted.

\case
In the remaining case $D''$ has at most $5$ clasps, but is disconnected
(that is, the connected sum of two genus 1 diagrams).
\end{caselist}

In all 3 cases the classification of genus two
diagrams of \cite{gen2} shows that $D''$ has at most 13 crossings.
Thus $D'$ has at most 14 crossings, and $D$ has at most 16 crossings.

Now, the computer check of all $\le 16$ crossing diagrams without
(reverse or parallel) clasp (up to flypes they are 203) shows
that the ratio $lk/v_2$ on genus 3 diagrams is indeed maximal on
the mentioned 9 crossing diagram of $8_{19}$. \qed

Now we use lemma \reference{LM2} to obtain an inequality between
$v_2$ and $c$. To make the expressions somewhat simpler,
instead of $3\mybrtwo{c-1}$ in the following we use $\myfrac{4}{3}c$,
which is weaker except for $c\le 8$ or $c=10,12,14,16$. For these
crossing numbers $lk\ge \myfrac{4}{3}c$ is checked directly by computer
and is found to be true except for $c=3,4$, in which cases the
subsequent claims are matter of straightforward verification.

Therefore, from lemma \ref{LM2} we obtain the following 

\begin{corr}
For any positive knot $K$ of genus at most $3$, we have
$v_2(K)/c(K)\ge 20/63$. \qed
\end{corr}

To get from the positive to the almost positive case, we have
the following lemma.

\begin{lem} \label{LM1}
If $D'$ is an almost positive diagram that by flypes cannot be transformed
into a diagram with (the negative crossing involved into) a resolved clasp,
and $D$ is obtained from $D'$ by switching the negative crossing to a
positive one, then $v_2(D')\ge v_2(D)-\myfrac{1}{5}lk(D)$.
\end{lem}

\proof Let $p$ be the negative crossing in $D'$. First we show that
if $a\cap p$, then $\#\{\,q\,:q\cap a\,\}\ge 4$.

Assume the contrary, i.~e., $\ex a\cap p\,:\,\ex! c\ne p\,:\,c\cap a$.
If now $\ex d\cap c:\,d\ncap p$, then $2C(d,c,a)$ would give a third
chord intersecting $a$, and in the same way $2C(d,p,a)$ will do if
$\ex d\cap p:\,d\ncap c$.
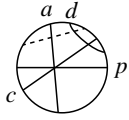
\[
\CD{
  \labch{180}{0}{p}
  \labch{265}{105}{a}
  \labch{40}{210}{c}
  \piccurve{80 1 x polar}{80 1 x polar 0.2 -}
	   {20 1 x polar x 0.2 - x}{20 1 x polar}
  \picputtext{80 1.3 x polar}{$d$}
  {\piclinedash{0.1}{0.05}
   \chrd{60}{150}
  }
}
\]
Therefore, $\{\,q\,:q\cap c\,\}=\{\,q\,:q\cap p\}$, which means that
after possible flypes $p$ and $c$ form a (resolved) clasp.

Therefore, $\fa a\cap p\,:\,\#\{\,q\,:q\cap a\,\}\ge 4$. This implies that,
if we set $v(p):=\#\{\,a\,:a\cap p\,\}$, then there are at least $\myfrac{3}
{2}v(p)$ linked pairs in $D'$ not involving $p$. The claim then follows from the
fact (which is straightforward from the Gau\ss{} sum formula), that if
$D$ is the (positive) diagram obtained from $D'$ by switching $p$, then
$v_2(D)-v_2(D')=\myfrac{v(p)}{2}$. \qed

\begin{corr}
For any almost positive knot $K$ of genus 2, we have 
\[
\frac{v_2(K)}{c(K)}\,\ge \,\frac {16}{315}\,.
\]
\end{corr}

\proof An almost positive diagram $D$ of $K$ has genus at most 3.
Thus using proposition \ref{PP1} and lemma \ref{LM1}, we get
\[
\frac{v_2(D)}{lk(D)}\,\ge\,\frac 5{21}-\frac 15\,=\,\frac {4}{105}\,.
\]
Multiplying the r.h.s. by $\myfrac{4}{3}$ and replacing $lk(D)$
by $c(D)$ (by the weaker version of lemma \reference{LM2}),
and then $c(D)$ by $c(K)$, gives the result. \qed

Note, that we lose the genus three case, as there may be almost
positive genus three knots that only possess almost positive diagrams
of genus four. On the other hand, the satement for the genus one case
is obsolete, because almost positive genus one knots do not exist
(see \cite{gen2}).

A presumably difficult combinatorial question that makes sense to
ask now is:

\begin{ques}\label{qq}
Does $lk/v_2$ grow unboundedly on positive diagrams?
\end{ques}

\begin{figure}[htb]
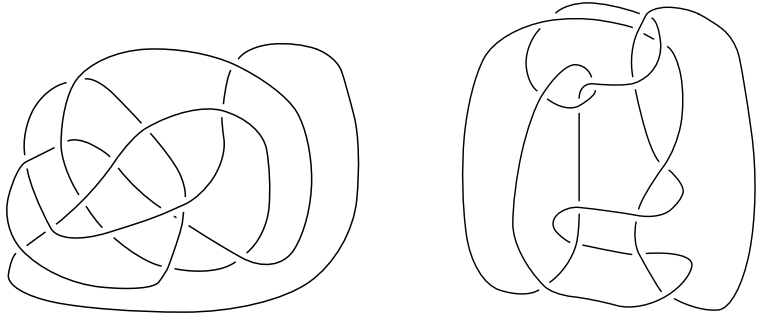

\[
\epsfs{4cm}{t-15-126448-1} \hskip1.3cm
\epsfs{4cm}{t-15-126448}
\]
\caption{\label{figcex}The 16 crossing diagram of the knot $15_{126448}
$, maximizing the ratio $lk/v_2$ on all positive diagrams of at most 16
crossings, and the diagram of the same knot included in Thistlethwaite's
table.}
\end{figure}

\begin{rem}
The computer experiment revealed that the maximal value of $lk/v_2$
on positive $\le 16$ crossing diagrams is attained on the diagram on
figure \ref{figcex} of the (4,5)-torus knot $15_{126448}$, where it is
$\myfrac{64}{15}=4.2\ol{6}$. This, of course, is far from giving any
indication about a (positive or negative) answer for question \ref{qq}.
\end{rem}

Coming back to the case of arbitrary genus, we can obtain some partial
improvements of results in \cite{pos} and such proved above in this
paper.

\begin{corr}\label{CoR1}
For any $\dl<\myfrac{3}{8}$ and any $g\in \bN$, there are only
finitely many positive knots of genus $g$ with $v_2<\dl c$.
\end{corr}

\proof Use lemma \ref{LM2} and theorem \ref{T3}. \qed

\begin{corr}\label{COR1}
For any $\dl<\myfrac{3}{40}$ and any $g\in \bN$, there are only
finitely many (almost positive diagrams of) almost positive knots 
of genus $g$ with $v_2<\dl c$.
\end{corr}

\proof 
By the extension of Bennequin's inequality \cite[theorem 3]{Bennequin}
to arbitrary diagrams (see \cite{Rudolph,pos}), an almost positive
knot of genus $g$ has almost positive diagrams only of genus $g$
or $g+1$, and thus it suffices
to consider diagrams of bounded genus. It this situation 
theorem \ref{T3} shows that, up to finitely many exceptions, on diagrams
of given genus $v_2>(\myfrac{1}{4}-\eps)lk$. Therefore, again
up to finitely many exceptions, on almost positive diagrams of given
genus, $v_2>(\myfrac{1}{4}-\myfrac{1}{5}-\eps)lk$ by lemma \ref{LM1}.

The above constant $\myfrac{3}{40}$ comes from the multiplication
of $\myfrac{1}{4}-\myfrac{1}{5}=\myfrac{1}{20}$ by $\myfrac{3}{2}$,
occurring in lemma \reference{LM2} (the defect constant
in $c$ therein is put into the $\eps$). Applying lemma
\ref{LM2} we need to remark that if the positive diagram admits a move 
of the type \eqref{second}, then so it does after the crossing change,
unless then it looks like
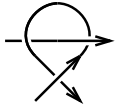
\[
\diag{1cm}{2}{1.6}{
  \picline{0.3 1}{1 1}
  \picmulticurve{0.12 1 -1.0 0}{1.0 1.5}{0.7 1.5}{0.4 1.1}{0.7 0.8}
  \picvecline{0.7 0.8}{1.3 0.2}
  \picmultivecline{0.12 1 -1.0 0}{0.7 0.2}{1.3 0.8}
  \piccurve{1.0 1.5}{0.7 1.5}{0.4 1.1}{0.7 0.8}
  \picmulticurve{0.12 1 -1.0 0}{1.3 0.8}{1.6 1.1}{1.3 1.5}{1.0 1.5}
  \picmultivecline{0.12 1 -1.0 0}{1 1}{1.7 1}
}\,.
\]
But in this case the negative crossing is linked with only two
positive crossings, and therefore $v_2$ decreases only by one under
the crossing change.  \qed

If we use the corollary for diagrams, we obtain

\begin{corr}
Each almost positive knot has only finitely many almost positive
diagrams.
\end{corr}

\proof If there were infinitely many such diagrams, they would be
of bounded genus, and then by corollary \reference{COR1}, we would
have $v_2\to\infty$ on them, a contradiction. \qed

The next consequence is the main step towards the proof of
theorem \reference{tL}.

\begin{corr}\label{cPV}
There are only finitely many almost positive knots with
given $\max\deg V$ or $\max\deg_v P$ (in the variables of
\cite{MorCro}). That is, if $K_1,K_2,\dots$ is a sequence of distinct
almost positive knots, then $\max\deg V_{K_i}\to\infty$.
\end{corr} 

\proof As for almost positive knots by \cite[corollary 4.4]{restr}
we have $\min\deg V\ge g-1$, assume, taking a subsequence,
that the genera of $K_i$ are bounded. Then by \cite[corollary 7.1]{beha}
any coefficient of $V$ admits only finitely many values on $\{K_i\}$,
and if $\max\deg V$ were bounded on $\{K_i\}$, so would be $v_2=-\myfrac
{1}{6}V''(1)$, contradicting corollary \ref{COR1}.

This argument establishes the desired property for a subsequence of
$(K_i)$, but applying it to any subsequence of $(K_i)$ gives the result.

Exactly the same reasoning applies for $P$ (in fact, even for its
absolute $[P]_{z^0}$ or quadratic terms $[P]_{z^2}$ of the Alexander
variable, which contain $v_2$, see \cite{KanMiy}). \qed

Recall that the Thurston--Bennequin number $tb(\cK)$ of a Legendrian
knot $\cK$ in the standard contact space $(\bR^3(x,y,z),\,dx+y\,dz)$
is the linking number of $\cK$ with $\cK'$, where $\cK'$ is obtained
from $\cK$ by a push-forward along a vector field transverse to the
(hyperplanes of the) contact structure.

The Maslov (rotation) index $\mu(\cK)$ of $\cK$ is the degree of the map
\[
t\in S^1\,\mapsto\,\frac{\pr\,\frac{\partial \cK}{\partial t}(t)}{
\bigl |\pr\,\frac{\partial \cK}{\partial t}(t)\bigr |}\in S^1\,,
\]
where $\pr\,:\,\bR^3\to\bR^2\simeq \bC$ is the projection
$(x,y,z)\mapsto(y,z)$.

\proof[of theorem \reference{tL}] Combine the result of corollary
\reference{cPV} for $P$ with the inequalities of \cite[\S 2, theorem
2.4]{TabFuchs} for $tb$ and $\mu$ coming from $\min\deg_vP$ (take care
of the mirroring convention for the contact structure). \qed

\begin{rem}
The estimate of the degrees of $V$ and $P$ can be made, if desired
(in particular with regard to theorem \reference{tL}), more explicit.
The lower bound we obtain is roughly of the form $g(K)+\sqrt[C_{g(K)}]
{c(K)}$ (for $K$ positive or almost positive), with certain constants
$C_g$ depending on $g$ only, which will need to incorporate the
numbers $d_g$ of \cite{gen1} and the exceptional knots of
corollaries \reference{CoR1} and \reference{COR1}. See also \cite{gwg},
where the estimate for $V$ and for positive knots was given.
\end{rem}

\begin{corr}
For any given genus $g$, there are only finitely many almost positive
knots of genus $g$ with given $\Dl$ or $Q$ polynomial.
\end{corr} 

\proof Use again that $v_2$ is contained in $\Dl$ and $Q$. \qed

\begin{rem}
Of course in view of the results on positive knots it is reasonable
to conjecture the statement to be true even without the genus
condition. For the Alexander polynomial it would follow from the
conjecture made in \cite{restr}:
\end{rem}

\begin{conj}
If a knot $K$ is almost positive, then $\max\deg\Dl_K=\min\deg V_K$.
\end{conj}

In \cite{pos} we used the inequality of theorem 6.1 therein to show that
if a positive knot has a (reduced) positive diagram of $c$ crossings,
then its crossing number is at least $\sqrt{2c}$. This result can now be
extended to almost positive knots, if we fix the genus.

\begin{corr}\label{COR3}
Let $\dl<\myfrac{3}{5}$. Then for any genus $g$ there are only finitely
many almost positive knots $K$ with a (reduced) almost positive diagram
of $c$ crossings but $c(K)<\sqrt{\dl c}$.
\end{corr}

\proof Use again (the proof of) corollary \ref{COR1} and theorem
2.2.E of \cite{VirPol2}. The constant $\myfrac{3}{5}$ comes from
multiplying the $\myfrac{1}{20}$ appearing in the proof of corollary
\ref{COR1} by the $8$ in the denominator of theorem 2.2.E of
\cite{VirPol2}, giving $\myfrac{2}{5}$,
and using the inequality $lk(D)\ge \bigl(\myfrac{3}{2}-\eps\bigr)
c(D)$ for a reduced diagram $D$. \qed

This result is indeed a little technical, but is exactly the missing
piece required to extend the results of \cite{gen1} that the number
of positive knots of given genus or unknotting number grows polynomially
in the crossing number.

\proof[of theorem \reference{tgu}]
This is almost a repetition of the proof in \cite{gen1}. Consider
first the genus case and fix $g$. Then, by Bennequin's inequality,
we need to consider diagrams of genus at most $g+1$. Then  theorem 3.1
of \cite{gen1} shows that the number of such diagrams grows polynomially
in the crossing number $c(D)$ \em{of the diagram}. But now, as we
consider diagrams of bounded genus, we have by corollary \ref{COR3}
only finitely many exceptions to throw out, in order to be allowed to
apply, say, the inequality $c(K)\ge \sqrt{\myfrac{1}{2}c(D)}$, for $D$ a
diagram of $K$. This shows, that (up to these finitely many exceptions)
we obtain all knots $K$ of given $c(K)$ by taking diagrams of at most
$c(D)\le 2c(K)^2$ crossings, and a polynomial in $2c(K)^2$ is a
polynomial (of double degree) in $c(K)$. 

The unknotting number result follows from that for the genus, because
$u\ge g-1$ for an almost positive knot. \qed

A final application of the methods described in this section is
a simple proof of the following fact.

\begin{prop}
There are only finitely many almost
positive knots $K$ of given genus $g(K)$ and given braid index $b(K)$.
\end{prop}

\proof Because of the braid index inequality 
\[
b(K)\ge 1+(\max\deg_v P_K-\min\deg_v P_K)/2
\]
of Franks--Williams \cite{WilFr} and Morton \cite{Morton} and
corollary \reference{cPV}, it suffices to show that for
almost positive knots $K$ of given genus $g$, $\min\deg_v P_K$
is bounded (above). Now, it follows from the identity 
$P_K(v,v+v^{-1})=1$ that $\min\deg_v P_K\le \max\deg_m P_K$,
from an inequality of \cite{Morton} that $ \max\deg_m P_K\le 2\tl g(K)$,
and from Bennequin's inequality that $\tl g(K)\le g(K)+1$ for an
almost positive knot $K$. Joining all this, we obtain
\[
\min\deg_v P_K\le \max\deg_m P_K\le 2\tl g(K)\le 2g(K)+2\,,
\]
and thus the desired bound. \qed


%
%

\section{The signature\label{Ssig}}

Although the Vassiliev invariant inequalities for positive knots
of \cite{pos} have nice theoretical consequences, they reveal
too weak to exclude knots with almost positive diagrams from
being positive. In this regard the most handy criterium is due to
Przytycki-Taniyama \cite{PrzTan} using the signature $\sg$. We apply
this criterium to show that there are infinitely many almost positive
knots. Very unfortunately, their work is unpublished, and so we give 
an independent proof of their results, using some methods and results of
\cite{pos}, \cite{gen1} and \cite{gen2}.

\begin{theorem} (see \cite{PrzTan})
Any positive knot has signature at least 4, except for the 
$(p,q,r)$-pretzel knots with $p,q,r>0$ odd, which have signature 2.
\end{theorem}

\proof In \cite[\S 6]{pos} we introduced a move on (positive)
diagrams, called \em{loop-move} to show an inequality of the Casson
invariant (theorem 6.1) and observed that it can also be used to show 
the positivity of the signature, a result of Cochran and Gompf
\cite[corollary 3.4, p.~497]{CochranGompf} and Traczyk \cite{Traczyk}.

The loop-move from a diagram $A$ to a diagram $B$, henceforth
denoted by $A\to B$, consists in choosing a segment of the line in $A$
between the two passings of a crossing, such that it has no
self-crossings, and removing of this segment by switching half of the
crossings on it (and elimination of all reducible crossings thereafter).
Such a move never augments the signature, that is, $\sg(B)\le\sg(A)$,
and a finite number of such moves makes any positive diagram trivial,
that is, for any diagram $D$, there is a sequence $D=D_0\to D_1\to
\dots\to D_n$ of loop moves, with $D_n$ being the zero crossing diagram.

By considering a loop move sequence in which $n$ is maximal,
the diagram $D':=D_{n-1}$ before applying the last step of such a
sequence can be chosen so that \em{any} loop move unknots it.
If $\sg(D')\ge 4$, then the assertion is satisfied. By \cite[exercise
6.4]{pos} and \cite{gen1} $\sg(D')=2$ exactly if $D'$ is a diagram of
genus one, in which case it corresponds to one of the desired pretzel
knots.

Therefore, it remains to show that if $D'$ is of genus one, but
$n>1$, then $\sg(D)\ge 4$. We show that already $\sg(D'')\ge 4$ for
$D'':=D_{n-2}$.

\begin{figure}[htb]
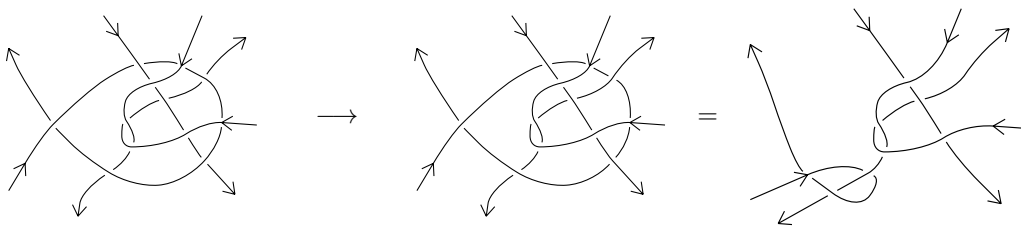

\[
\epsfsv{3.0cm}{t-postang}\qquad\lra\qquad
\epsfsv{3.0cm}{t-postang2}\quad=\quad
\epsfsv{3.0cm}{t-postang3}
\]
\caption{\label{figaloop}An almost loop-move. Take a crossing whose
smoothing out gives a component with no self-crossings, and switch
some of the crossings on the segment except the last one, so that 
the segment to (over-/un\-der-)\-pass them all in different way
than the last one. Then the segment can be simplified to have just
two cossings on it.}
\end{figure}

Replace the loop move $D''\to D'$ by an almost loop-move, as shown
on figure \ref{figaloop}. One obtains from $D''$ a (positive) diagram
$\tl D$, from which $D'$ arises by switching on of the two crossings
of the segment removed by $D''\to D'$. By direct observation one sees
that (independently of the orientation), the move $\tl D\to D'$
(which is formally also a loop-move), preserves the number of Seifert
circles, and hence $\tl D$ has genus two.

But all (\em{inter alia}, positive) genus two diagrams
are classified in \cite{gen2}, and are in particular shown to be
transformable by changing positive crossings to certain 24
(called therein ``generating'') diagrams, and it is straightforward
to check that for all them $\sg\ge 4$ (see corollary 3.2 therein),
which shows the assertion. \qed

\begin{rem}
When the length $n$ of the loop move sequence becomes large, then
it is apparent that a much better estimate should be obtainable for the
signature. This evidence is compatible with the conjecture made in
\cite{2apos} that $\sg$ is bounded below by an increasing function of
the genus for positive (and hence also for almost positive) knots.
To prove such a statement using loop moves only seems difficult, though.
Although a loop move in general reduces $\sg$, this needs not always
be the case. Worse yet, there are positive diagrams, on which \em{no
one} of the applicable loop moves strictly reduces $\sg$ (one such
diagram is the 14 crossing diagram of the knot $14_{45657}$ given in
\cite{2apos}).
\end{rem}

\begin{figure}[htb]
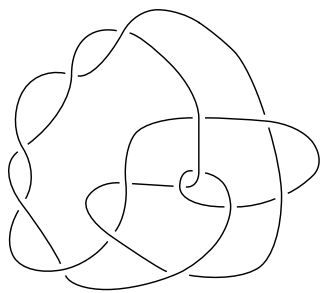

\[
\epsfs{4cm}{k-14-34605}
\]
\caption{\label{fig14-34605}The knot $!14_{34605}$, a member of an
infinite family of almost positive knots.}
\end{figure}

\begin{exam}
There are infinitely many almost positive knots.
Consider the knot $!14_{34605}$ of Thistlethwaite \cite{HTW,KnotScape}
on figure \ref{fig14-34605}. It is a genus two knot of determinant $-3$,
and hence has signature $2$. Switching one of the five crossings
forming reverse clasps in the upper left part of the diagram, we obtain
$!12_{1692}$, a knot with the same genus and determinant, and also
the same for $!10_{145}$, obtained by switching two of these five
crossings. Therefore, applying $\bt$ moves at the five crossings
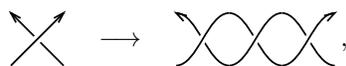
\[
\diag{7mm}{1}{1}{
    \picmultivecline{-5.0 1 -1.0 0}{1 0}{0 1}
    \picmultivecline{-5.0 1 -1.0 0}{0 0}{1 1}
}\quad\lra\quad
\diag{7mm}{3.0}{2}{
  \picPSgraphics{0 setlinecap}
  \pictranslate{0.99 1}{
    \picrotate{-90}{
      \rbraid{0 -0.5}{1 1}
      \rbraid{0 0.5}{1 1}
      \rbraid{0 1.5}{1 1}
      \pictranslate{-0.5 0}{
      \picvecline{0.03 1.95}{0 2}
      \picvecline{0.03 -.95}{0 -1}
    }
    }
    }
}\,,
\]
we obtain an infinite family of knots starting with $!10_{145},$
$!12_{1692},$ $!14_{34605},$ $!16_{970714},\dots$ with genus two and
determinant $-3$, and hence signature $2$. If some of them were
positive, then it
would have to be a pretzel knot, and would have genus one. On the other
hand, the diagrams are evidently all almost positive, therefore
so are all these knots. Moreover, the diagrams obtained are indeed of
minimal crossing number as can be shown by examining their $Q$
polynomial \cite{Kidwell} and using the fact that by \cite[\S 5]{restr}
the knots are non-alternating.
\end{exam}

\begin{corr}
Any almost positive knot has positive signature.
\end{corr}

\proof Switch the one negative crossing to the positive, which
does not augment $\sg$ by more than two, and apply the previous theorem.
The only point to remark is that an almost positive
diagram of genus one corresponds to a positive knot. \qed

\begin{rem}
In fact, it is interesting to remark that the loop move can also be
used to show Taniyama's theorems of \cite{Taniyama}, that any
non-trivial
knot diagram can be crossing switched to a diagram of the trefoil,
and that any connected knot diagram of genus at least 2 can be done so
to the (5,2)-torus knot. For this we just need not to pull out the
loops in the loop moves and to switch any subsequent loop on top of this
additional set of crossings (remaining by not pulling out the previous 
loops). Then we are left with the 24 generating genus two diagrams
of \cite[fig. 5]{gen2} (modulo flypes, which are irrelevant here),
for which the claims can be checked directly.
\end{rem}

\section{A conjecture}

The considerations open a natural question about further
generalizations.

\begin{defi}
A knot is called $2$-almost positive, if the minimal number of negative
crossings in all its diagrams is $2$.
\end{defi}

Although it is not clear how the classification result
should carry over to $2$-almost positive unknot diagrams, at first
glance it appears that at least theorem \reference{th1}
should generalize to this case. Surprisingly, this turns out
not to be the case.

\begin{exam}
The knots $!6_1$ and $!6_2$ have in their ($2$-almost positive) 6
crossing 
diagrams $v_3=-4$. This also shows that $!6_1$ and $!6_2$ are indeed
$2$-almost positive. Furthermore, as $!6_1$ and $!6_2$ can be unknotted
in their 6 crossing diagrams by switching only positive crossings, this
shows, that (although measuring positivity in general by positive
values, contrarily to the signature) $v_3$ increases sometimes, when
a positive crossing is switched to a negative one.
\end{exam}

Therefore, unfortunately, our approach will very unlikely carry over to
classify $2$-almost positive unknot diagrams. This has been achieved
in \cite{gen2} using a deeper tool~-- the version of the inequality of
Bennequin for arbitrary knot diagrams \cite{Rudolph,pos}.

However, with some heuristics the above example leads
to the following conjecture:

\begin{conj}
Let $L$ be a $2$-almost positive even crossing number diagram
minimizing $v_3$ over all diagrams of that crossing number.
Then $L$ is a diagram of a $(a_1,\dots,a_k,2)$ pretzel knot,
$a_i\in\{1,3\}$.
\end{conj}

\noindent{\bf Acknowledgements.}
I would like to thank to V.~Chernov for some helpful remarks and
discussions and for his invitation to Zurich in summer 1999.

\end{document}

%% file: myeqn.tex
\newenvironment{myeqn*}[1]{\begingroup\def\@eqnnum{\reset@font\rm#1}%
\xdef\@tempk{\arabic{equation}}\begin{equation}\edef\@currentlabel{#1}}
{\end{equation}\endgroup\setcounter{equation}{\@tempk}\ignorespaces}

\newenvironment{myeqn}[1]{\begingroup\let\eq@num\@eqnnum
\def\@eqnnum{\bgroup\let\r@fn\normalcolor 
\def\normalcolor####1(####2){\r@fn####1#1}%
\eq@num\egroup}%
\xdef\@tempk{\arabic{equation}}\begin{equation}\edef\@currentlabel{#1}}
{\end{equation}\endgroup\setcounter{equation}{\@tempk}\ignorespaces}